\documentclass[10pt]{amsart}
\usepackage[T1]{fontenc}
\usepackage{a4wide}
\usepackage[latin1]{inputenc}
\usepackage{amsmath}
\usepackage{amsfonts}
\usepackage{amssymb}
\usepackage{amsthm}
\usepackage{subfigure}
\newtheorem{theo}{Theorem}[section]
\newtheorem{prop}[theo]{Proposition}
\newtheorem{coro}[theo]{Corollary}
\newtheorem{lemm}[theo]{Lemma}

\theoremstyle{definition}
\newtheorem{def1}[theo]{Definition}
\theoremstyle{remark}
\newtheorem*{rema}{Remark}
\newcommand{\Op}{\text{Op}}

\title[Entropy of quantum limits for symplectic linear maps of $\mathbb{T}^{2d}$]{Entropy of quantum limits for symplectic linear maps of the multidimensional torus}
\author{Gabriel Rivière}
\address{Centre de Math\'ematiques Laurent Schwartz (UMR 7640), \'Ecole Polytechnique, 91128 Palaiseau Cedex, France}
\email{gabriel.riviere@math.polytechnique.fr}

\begin{document}

\begin{abstract}
In the case of a linear symplectic map $A$ of the $2d$-torus, we give a lower bound for the entropy of semiclassical measures. In this case, semiclassical measures are $A$-invariant probability measures associated to sequences of high energy quantum states. Our main result is that, for any \emph{quantizable} matrix $A$ in $Sp(2d,\mathbb{Z})$ and any semiclassical measure $\mu$ associated to it, one has
$$h_{KS}(\mu,A)\geq \sum_{\beta\in\text{sp}(A)}\max\left(\log|\beta|-\frac{\lambda_{\max}}{2},0\right),$$
where the sum is taken over the spectrum of $A$ (counted with multiplicities) and $\lambda_{\max}$ is the supremum of $\{\log |\beta|:\ \beta\in\text{sp}(A)\}$. In particular, our result implies that if $A$ has an eigenvalue outside the unit circle, then a semiclassical measure cannot be carried by a closed orbit of $A$.
\end{abstract}

\maketitle

\section{Introduction}

The semiclassical principle asserts that in the high energy limit of quantum mechanics, phenomena from classical mechanics appear. One of the issue in quantum chaos is to understand the influence of the chaotic properties of the dynamical system (Anosov property, ergodicity, etc.) on the asymptotic behaviour of the eigenvectors of the quantum propagator.\\
The first result in this direction is due to Shnirelman~\cite{Sc}, Zelditch~\cite{Ze} and Colin de Verdière~\cite{CdV}. It states that, given an ergodic geodesic flow on a Riemannian manifold $M$, \emph{almost all} eigenfunctions of $\Delta$ are equidistributed on $S^*M$ in the high energy limit. This phenomenon is known as quantum ergodicity and has many extensions. Rudnick and Sarnak formulated the so-called `Quantum Unique Ergodicity Conjecture' which states that for manifolds of negative curvature, all the eigenfunctions of $\Delta$ are equidistributed in the high energy limit~\cite{RS}. This conjecture remains widely open in the general case.\\
To study quantum chaos, an approach is to study toy models, i.e. simple symplectic dynamical systems that are highly chaotic and that can be quantized in a standard way. One of the main advantage of such a system is that we can make explicit computations. Among the several toy models is the family of symplectic linear automorphisms on the $2d$-torus $\mathbb{T}^{2d}$. In fact, a matrix $A$ in $Sp(2d,\mathbb{Z})$ which does not have $1$ as an eigenvalue can be quantized in a standard way~\cite{BodB2}. We will say that $A$ is \emph{quantizable}. In this case, the phase space $\mathbb{T}^{2d}$ is compact and in particular, the natural Hilbert space that can be associated to the level $N$ will be finite dimensional. It will be denoted $\mathcal{H}_N(\kappa)$ (see section~\ref{QM-torus}) and will be of dimension $N^d$. The semiclassical parameter is denoted $\hbar$ and satisfies $2\pi\hbar N=1$ (where $N$ is an integer). The set of classical observables will be the set of smooth functions on the torus $\mathcal{C}^{\infty}(\mathbb{T}^{2d})$. There exists then a positive quantization procedure $\Op_{\kappa}^{AW}(.)$ that associates to each observable $a$ a linear operator $\Op_{\kappa}^{AW}(a)$ on $\mathcal{H}_N(\kappa)$. This procedure is called the anti-Wick quantization and is constructed from the family of coherent states~\cite{BdB}. Moreover, there is a quantum propagator $M_{\kappa}(A)$ corresponding to $A$ which acts on $\mathcal{H}_{N}(\kappa)$. This propagator satisfies the Egorov property:
\begin{equation}\label{Egorov-torus}M_{\kappa}(A)^{-1}\Op_{\kappa}^{AW}(a)M_{\kappa}(A)=\Op_{\kappa}^{AW}(a\circ A)+O_a(N^{-1}).\end{equation}
For any eigenvector $\varphi^N$ of $M_{\kappa}(A)$ in the energy level $\mathcal{H}_{N}(\kappa)$, one can define the following measure on the torus:
\begin{equation}\label{def-measures-torus}
\tilde{\mu}^N(a):=\langle\varphi^N|\Op_{\kappa}^{AW}(a)|\varphi^N\rangle_{\mathcal{H}_N(\kappa)}.
\end{equation}
This quantity gives a description of the quantum state at energy $N$ in terms of the position and the velocity, i.e. of the quantum state and its $N$-Fourier transform. Thanks to Egorov theorem, we have that any weak limit of the corresponding $(\tilde{\mu}^N)$ on the high energy limit (i.e. as $N$ tends to infinity) is an $A$-invariant measure on the torus. We call these accumulation points semiclassical measures. In this setting, Bouzouina and de Bièvre proved an analogue of Shnirelman's theorem~\cite{BdB}, i.e. \emph{almost all} the $\tilde{\mu}^N$ converge weakly to the Lebesgue measure on the torus if $A$ is \emph{ergodic}\footnote{It means that no eigenvalue of $A$ is a root of unity.}. However, it has been shown that the Quantum Unique Ergodicity property fails in this setting and in any dimension~\cite{FNdB},~\cite{Ke}. Precisely, given an \emph{hyperbolic}\footnote{It means that $A$ has no eigenvalue on the unit circle.} matrix $A$, Lebesgue measure is not the only accumulation point of the $\tilde{\mu}^N$. For $d=1$, it was proved by de Bièvre, Faure and Nonnenmacher that $\frac{1}{2}(\delta_0+\text{Leb})$ is a semiclassical measure~\cite{FNdB}. In higher dimensions and under arithmetic assumptions on $A$, Kelmer constructed semiclassical measures supported on submanifolds of $\mathbb{T}^{2d}$~\cite{Ke}.\\
Even if we know that the set of semiclassical measures is not reduced to the Lebesgue measure for quantized maps of the torus, one can ask about the properties of these semiclassical measures. For instance, it was shown in~\cite{BodB2} and~\cite{FN} that if we split the semiclassical measure into its pure point, Lebesgue and singular continuous components, $\mu=\mu_{\text{pp}}+\mu_{\text{Leb}}+\mu_{\text{sc}}$, then $\mu_{\text{pp}}(\mathbb{T}^2)\leq\mu_{\text{Leb}}(\mathbb{T}^2)$ and in particular $\mu_{\text{pp}}(\mathbb{T}^2)\leq 1/2$.

\subsection{Statement of the main theorem}

In~\cite{An}, Anantharaman proved that for a compact riemannian manifold $M$ of Anosov type, the Kolmogorov-Sinai entropy of any semiclassical measure associated to a sequence of eigenfunctions of $\Delta$ is positive (see section~\ref{s:proof-maintheo} or~\cite{Wa} (chapter $4$) for a definition of the entropy). Her result proves in particular that eigenfunctions of the Laplacian cannot concentrate only on a closed
geodesic, in the large eigenvalue limit. Translated in the context of our toy models, her result says that, for any symplectic and \emph{hyperbolic} matrix $A$, the Kolmogorov-Sinai entropy of a corresponding semiclassical measure is positive. In particular, a semiclassical measure cannot be supported only on a closed orbit of $A$. In subsequent works with Koch and Nonnenmacher~\cite{AN2},~\cite{AKN}, they gave quantitative lower bounds on the Kolmogorov-Sinai entropy of semiclassical measures. In the model of quantized maps of the $2d$-torus, their result can be written, for any semiclassical measure associated to an \emph{hyperbolic} matrix $A$:
\begin{equation}\label{e:AN-result}
h_{KS}(\mu,A)\geq\sum_{i=1}^{2d}\max\left(\log|\beta_i|,0\right)-\frac{d}{2}\lambda_{\max},
\end{equation}
where $\{\beta_{i}:1\leq i\leq 2d\}$ is the spectrum of $A$ (counted with multiplicities) and $\lambda_{\max}$ is the maximum of the $\log|\beta_i|$~\cite{AKN}. As $A$ is hyperbolic, it has exactly $d$ eigenvalues (counted with multiplicities) of modulus larger than $1$. We underline that their result was more general as it also deals with varying Lyapunov exponents.
One can remark that if $\lambda_{\max}$ is very large, the previous lower bound can be negative (and so the result empty). So they were lead to formulate that this bound should not be optimal. Regarding their result (and the different counterexamples), they conjectured that the optimal result would be that, for any semiclassical measure $\mu$ associated to an \emph{hyperbolic} matrix $A$,
\begin{equation}\label{e:conjecture}h_{KS}(\mu,A)\geq\frac{1}{2}\sum_{i=1}^{2d}\max\left(\log|\beta_i|,0\right).\end{equation}
Again, their conjecture was more general as they expected it to hold for situations where there are varying Lyapunov exponents. In this article, we will show:
\begin{theo}\label{t:main-theo} Let $A$ be a matrix in $Sp(2d,\mathbb{Z})$ such that $1$ is not an eigenvalue of $A$. Let $\mu$ be a semiclassical measure on $\mathbb{T}^{2d}$ associated to $A$. One has
\begin{equation}
h_{KS}(\mu,A)\geq\sum_{i=1}^{2d}\max\left(\log|\beta_i|-\frac{\lambda_{\max}}{2},0\right),
\end{equation}
where $\{\beta_i:1\leq i\leq 2d\}$ is the spectrum of $A$ and $\lambda_{\max}$ is the maximum of the $\log|\beta_i|$.
\end{theo}
A first comment is that we do not obtain exactly the lower bound expected by Anantharaman, Koch and Nonnenmacher. Compared with their result, the lower bound of our theorem improves their bound~\cite{AKN} and it always defines a nonnegative quantity. However, regarding the semiclassical measures constructed\footnote{These semiclassical measures have entropy equal to $\frac{1}{2}\sum_{i=1}^{2d}\max\left(\log|\beta_i|,0\right)$.} in~\cite{Ke}, the lower bound of our theorem should be suboptimal. We also recall that the Kolmogorov-Sinai entropy satisfies the Ruelle inequality, i.e. for any $A$-invariant measure $\mu$,
$$h_{KS}(\mu,A)\leq\sum_{i=1}^{2d}\max\left(\log|\beta_i|,0\right),$$
with equality if $\mu=\text{Leb}$~\cite{Wa} (chapter $8$). A second comment is that our theorem holds for any \emph{quantizable} matrix. Our assumption on the eigenvalues of $A$ is just made to ensure that the classical system can be quantized. We do not make any chaoticity assumption on the dynamical system $(\mathbb{T}^{2d},A)$ (for instance, we do not assume $A$ to be hyperbolic). In particular, our result implies that if $A$ has an eigenvalue outside the unit circle, then a semiclassical measure of $(\mathbb{T}^{2d},A)$ cannot only be carried by closed orbits of $A$.\\
\\
In the case of varying Lyapunov exponents, the lower bound~(\ref{e:conjecture}) has been shown to be true when one has only one Lyapunov exponent~\cite{G},~\cite{GR}. At this point, it is not clear to us how to combine both methods in order to obtain an explicit nonnegative lower bound on the entropy of semiclassical measures in a general setting. Finally, we underline that, in the case of hyperbolic automorphisms of $\mathbb{T}^2$, stronger results on the entropy of semiclassical measures were obtained by Brooks~\cite{Br} and that, in the case of locally symmetric spaces of rank $\geq2$, a similar lower bound was obtained by Anantharaman and Silberman~\cite{AS}.

\subsection{Strategy of the proof}

Compared with the original result of Anantharaman, the proof of inequality~(\ref{e:AN-result}) in~\cite{AKN} was simplified by the use of an entropic uncertainty principle due to Maassen and Uffink~\cite{MU}. This principle is a consequence of the Riesz-Thorin interpolation theorem and it can be stated as follows~\cite{AN2},~\cite{MU}:
\begin{theo}[Maassen-Uffink] Let $\mathcal{H}$ and $\tilde{\mathcal{H}}$ be two Hilbert spaces. Let $U$ be an unitary operator on $\tilde{\mathcal{H}}$. Suppose  $(\pi_i)_{i=1}^D$ is a family of operators from $\tilde{\mathcal{H}}$ to $\mathcal{H}$ that satisfies the following property of partition of identity:
$$\sum_{i=1}^D\pi_i^{\dagger}\pi_i=\text{Id}_{\tilde{\mathcal{H}}}.$$
Then, for any unit vector $\psi$, we have
\begin{equation}\label{ppe-incertitude-partition-chat} -\sum_{i=1}^D\|\pi_i\psi\|^2_{\mathcal{H}}\log\|\pi_i\psi\|^2_{\mathcal{H}}-\sum_{i=1}^D\|\pi_iU\psi\|^2_{\mathcal{H}}\log\|\pi_iU\psi\|^2_{\mathcal{H}}\geq-2\log\sup_{i,j}\|\pi_i U\pi_j^{\dagger}\|_{\mathcal{L}(\mathcal{H})}.
\end{equation}
\end{theo}
In~\cite{AKN}, the method was to use this principle for eigenfunctions of the Laplacian on $M$ and a well-chosen partition of $\text{Id}_{L^2(M)}$ so that the quantity in the left side of~(\ref{ppe-incertitude-partition-chat}) can be interpreted as the usual entropy from information theory~\cite{Wa}. One of the main difficulty (that already appeared in~\cite{An}) was then to give a sharp estimate on the quantity $\|\pi_i U\pi_j^{\dagger}\|_{L^2(M)\rightarrow L^2(M)}$ for this choice. In~\cite{GR3}, we managed to use the symmetries of the quantization procedure for quantized cat-maps on $\mathbb{T}^2$ to make a slightly different choice of partition. Precisely, we used the fact that the coherent states satisfy a property of partition of identity on $\mathcal{H}_N(\kappa)$ and we were able to implement this property when we applied the entropic uncertainty principle. With this remark, the quantity $\|\pi_i U\pi_j^{\dagger}\|_{\mathcal{H}_N(\kappa)\rightarrow \mathcal{H}_N(\kappa)}$ is easier to bound than the corresponding one in~\cite{AKN}. In fact, the bound can be derived from estimates on the propagation of coherent states under the quantum propagator as in~\cite{BodB} and~\cite{FNdB}.\\
\\
Our strategy will be to generalize our method for $d=1$ to higher dimensions. To do this, we would need to use a family of coherent states depending on the classical dynamic induced by $A$. In order to obtain the sharpest bound possible as in~\cite{GR3}, we would like to introduce coherent states which are not located in a ball but in an ellipsoid whose lengths will depend on the Lyapunov exponent of $A$ in every direction. It turns out that it is not possible to define such coherent states and to make the proof works. However, we will introduce a new quantization procedure that will be very similar to the anti-Wick procedure. We will construct this quantization using gaussian observables centered on ellipsoid adapted to the classical dynamic. We will prove that this procedure defines the same set of semiclassical measures and that it is adapted to the method developped in the case $d=1$. 

\subsection{Organization of the article}

In section~\ref{QM-torus}, we recall how the dynamical system $(\mathbb{T}^{2d},A)$ can be quantized in a standard way. In section~\ref{s:oseledets-sympl}, we collect some facts about the reduction of symplectic matrices. Then, in section~\ref{s:new-quant}, we construct a new quantization procedure adapted to the classical dynamic induced by $A$. In section~\ref{s:proof-maintheo}, we apply the entropic uncertainty principle to derive theorem~\ref{t:main-theo}. Finally, in section~\ref{main-estimate}, we prove a crucial estimate on our quantization procedure that we used in section~\ref{s:proof-maintheo}. This estimate is similar to the ones obtained for the propagation of coherent states in~\cite{BodB},~\cite{FNdB}. The appendices are devoted to the proof of crucial and technical lemmas that we admitted at different steps of the article.

\subsection*{Acknowledgements} I would like to thank warmly Nalini Anantharaman for introducing me to these kind of questions, for many discussions on this subject and for encouraging me to extend the method of~\cite{GR3} to the higher dimensional case. I also thank sincerely St\'ephane Nonnenmacher for his precious comments on a preliminary version of this paper. Finally, I would like to thank Fr\'ed\'eric Faure for several interesting discussions on the model of quantized cat-maps and Nicolas Vichery for helpful references on symplectic linear algebra. This work was partially supported by Agence Nationale de la Recherche under the grant ANR-09-JCJC-0099-01.

\section{Quantum mechanics on the $2d$-torus}
\label{QM-torus}
In this section, we recall some basic facts about quantization of linear symplectic toral automophisms. We follow the approach and notations used by de Bi\`evre $\&$ al. in previous articles~\cite{BodB2},~\cite{BdB}. We refer the reader to them for further details and references. We denote $\mathbb{T}^{2d}:=\mathbb{R}^{2d}/\mathbb{Z}^{2d}$ the $2d$-torus.

\subsection{Quantization of the phase space}
In physical words, $\mathbb{R}^{d}$ (or $\mathbb{T}^{d}$) is called the configuration space and $\mathbb{R}^{2d}$ (or $\mathbb{T}^{2d}$) is the phase space associated to it. In this article, $\rho=(x,\xi)$ will denote a point of the phase space, i.e. points of $\mathbb{R}^{2d}$ or $\mathbb{T}^{2d}$. The usual scalar product on $\mathbb{R}^{2d}$ is denoted $\langle a, b \rangle$ and $\sigma$ is the usual symplectic form on $\mathbb{R}^{2d}$, i.e. $\sigma(\rho,\rho')=\langle\rho,J\rho'\rangle$ where $\displaystyle J:=\left(\begin{array}{cc}0&-\text{Id}_{\mathbb{R}^d}\\ \text{Id}_{\mathbb{R}^d}&0\end{array}\right)$. For $\psi\in\mathcal{S}'(\mathbb{R}^d)$, we can define $(Q_j\psi)(x):=x_j\psi(x)$ and $(P_j\psi)(x):=\frac{\hbar}{\imath}\frac{\partial \psi}{\partial x_j}(x)$. This allows to define translation operators acting on a tempered distribution as:
$$U_{\hbar}(x,\xi):=e^{\frac{\imath}{\hbar}\sigma\left(\left(x,\xi\right),\left(Q,P\right)\right)}.$$
We underline that it is the standard representation of parameter $\hbar$ on $\mathcal{S}'(\mathbb{R}^d)$ of the Heisenberg group and that it is unitary on $L^2(\mathbb{R}^d,dx)$. Standard facts about this can be found in the book by Folland~\cite{Fo}. In particular, it can be shown that:
\begin{equation}\label{commutation-op-translation}U_{\hbar}(\rho)U_{\hbar}(\rho')=e^{\frac{\imath}{2\hbar}\sigma\left(\rho,\rho'\right)}U_{\hbar}(\rho+\rho').
\end{equation}
To define the quantum states associated to the phase space $\mathbb{T}^{2d}$, we can require that the Hilbert space has the same invariance under the translation operators $U_{\hbar}(q,p)$ for $(q,p)\in\mathbb{Z}^{2d}$. It means that the quantum states will have the same periodicity as the phase space. To do this, we let $\kappa=(\kappa_1,\kappa_2)$ be an element of $[0,2\pi[^{2d}$ and we require that, for all $(q,p)\in\mathbb{Z}^{2d}$, a quantum state $\psi$ should check the following condition:
$$U_{\hbar}(q,p)\psi=e^{\frac{\imath}{2\hbar}\langle q,p\rangle}e^{-\imath\langle\kappa_1,q\rangle+\imath\langle\kappa_2,p\rangle}\psi.$$
It can be remarked that $\kappa$ different from $0$ is allowed as, for $\alpha\in\mathbb{R}$, $\psi$ and $e^{\imath\alpha}\psi$ represent the same quantum state. The states in $\mathcal{S}'(\mathbb{R}^d)$ that satisfy the previous conditions are said to be the quantum states on the $2d$-torus and their set is denoted $\mathcal{H}_N(\kappa)$ where $N$ satisfies
\begin{equation}2\pi\hbar N=1.\end{equation}
It defines tempered distributions of period $1$ (modulo phase terms) whose $\hbar$-Fourier transform is also $1$-periodic. They are sums of Dirac distributions centered on a lattice of volume $\hbar^d$~\cite{BdB}. The following lemma can be shown~\cite{BdB}:
\begin{lemm}
$\mathcal{H}_N(\kappa)$ is not reduced to $0$ iff $N\in\mathbb{N}^*$. In this case, $\dim\mathcal{H}_N(\kappa)=N^d$. Moreover, for all $r\in\mathbb{Z}^{2d}$, $U_{\hbar}(\frac{r}{N})\mathcal{H}_N(\kappa)=\mathcal{H}_N(\kappa)$ and there is a unique Hilbert structure such that $U_{\hbar}(\frac{r}{N})$ is unitary for each $r\in\mathbb{Z}^{2d}$.
\end{lemm}
The Hilbert structure on $\mathcal{H}_N(\kappa)$ is not very explicit~\cite{BdB}. However, one can make it more clear using the following map which defines a surjection of $\mathcal{S}(\mathbb{R}^d)$ (Schwartz functions) onto $\mathcal{H}_N(\kappa)$:
\begin{equation}\label{definition-projecteur-chat}S(\kappa):=\sum_{(n,m)\in\mathbb{Z}^{2}}(-1)^{N\langle n,m\rangle}e^{\imath(\langle\kappa_1,n\rangle-\langle\kappa_2,m\rangle)}U_{\hbar}(n,m).\end{equation}
This projector associates to each state in $\mathcal{S}(\mathbb{R}^d)$ a state which is periodic in position and impulsion. Using it, we can define $|\phi,\kappa\rangle:=S(\kappa)|\phi\rangle$ and $|\phi',\kappa\rangle:=S(\kappa)|\phi'\rangle$ for $|\phi\rangle$ and $|\phi'\rangle$ in $\mathcal{S}(\mathbb{R}^d)$. Then, the following link between scalar products on $L^{2}(\mathbb{R}^d)$ and $\mathcal{H}_{N}(\kappa)$ holds:
\begin{equation}\label{produit-scalaire-tore}
\langle\kappa,\phi|\phi',\kappa\rangle_{\mathcal{H}_N(\kappa)}=\sum_{n,m\in\mathbb{Z}^{2}}(-1)^{N \langle n,m\rangle}e^{\imath(\langle\kappa_1,n\rangle-\langle\kappa_2,m\rangle)}\langle\phi|U_{\hbar}(n,m)|\phi'\rangle_{L^2(\mathbb{R})}.
\end{equation}
Finally, the following decomposition into irreducible subrepresentations of the discrete Weyl-Heisenberg group $\{(\frac{r}{N},\phi):r\in\mathbb{Z}^{2d},\phi\in\mathbb{R}\}$ can be written~\cite{BdB}:
$$L^2(\mathbb{R}^d)\cong\int_{[0,2\pi[^{2d}}\mathcal{H}_N(\kappa)d\kappa\ \text{and}\ U_{\hbar}\left(\frac{r}{N}\right)=\int_{[0,2\pi[^{2d}}U_{\kappa}\left(\frac{r}{N}\right)d\kappa.$$

\subsection{Weyl quantization}\label{ss:weyl-quant}
In the case of $\mathbb{R}^{2d}$, classical observables are functions of $\rho=(x,\xi)$ that belong to a certain class of symbols. We will use the following class of symbols:
$$S^k_{\nu}(1):=\{a\in\mathcal{C}^{\infty}(\mathbb{R}^{2d}):\ \text{for all multiindices}\ \alpha,\ \|\partial^{\alpha}a\|_{\infty}\leq \hbar^{-k-\nu|\alpha|} C_{\alpha}\},$$
where $\nu\leq\frac{1}{2}$. An usual way to quantize these observables is to use the Weyl quantization~\cite{DS},~\cite{EZ}. Let us recall the standard definition of this operator for an observable $a$:
$$[\Op_{\hbar}^{w}(a)u](x):=\frac{1}{(2\pi\hbar)^d}\int_{\mathbb{R}^{2d}}e^{\frac{\imath}{\hbar}\langle x-y,\xi\rangle}a\left(\frac{x+y}{2},\xi\right)u(y)dyd\xi.$$
We also recall that, from the Calder\'on-Vailancourt theorem (theorem $7.11$ in~\cite{DS} or theorem $4.22$ in~\cite{EZ}), we know that there exists an integer $D$ and a constant $C$ (depending only on $d$) such that
\begin{equation}\label{e:calderon-vail}\forall\ 0\leq\nu\leq\frac{1}{2},\ \forall a\in S^0_{\nu}(1),\ \left\|\Op_{\hbar}^w(a)\right\|\leq C\sum_{|\alpha|\leq D}\hbar^{\frac{|\alpha|}{2}}\|\partial^{\alpha}a\|_{\infty}.\end{equation}
We underline that the case $\nu=\frac{1}{2}$ is authorized for this last result~\cite{DS}. In the case of the $2d$-torus, classical observables are $\mathcal{C}^{\infty}$ functions on $\mathbb{T}^{2d}$ (and can be seen as a subset of $S^0(1)$). It can be shown that for $a\in\mathcal{C}^{\infty}(\mathbb{T}^{2d})$,
$$\Op_{\hbar}^{w}(a)=\sum_{r\in\mathbb{Z}^{2d}}a_{r}U_{\hbar}\left(\frac{r}{N}\right),$$
where $a_{r}$ is the $r$ coefficient of the Fourier serie of $a$, i.e. $a(\rho)=\sum_{a\in\mathbb{Z}^{2d}}a_{r}e^{-2\imath\pi\langle Jr,\rho\rangle}.$
Using the fact that $U_{\hbar}(r)\Op_{\hbar}^{w}(a)U_{\hbar}(r)^*=\Op_{\hbar}^{w}(a)$ (thanks to~(\ref{commutation-op-translation})), it follows that $\Op_{\hbar}^{w}(a)\mathcal{H}_N(\kappa)\subset\mathcal{H}_N(\kappa).$
In view of this remark, we shall denote $\Op_{\kappa}^{w}(a)$ the restriction of $\Op_{\hbar}^{w}(a)$ to $\mathcal{H}_N(\kappa)$. Finally, the following decomposition holds:
$$\Op_{\hbar}^{w}(a)=\int_{[0,2\pi[^{2d}}\Op_{\kappa}^{w}(a)d\kappa.$$
Recall from~\cite{ReSi} (theorem XIII.83) that such a decomposition implies:
\begin{equation}\label{e:relation-norm-op}\sup_{\kappa}\|\Op_{\kappa}^{w}(a)\|_{\mathcal{L}(\mathcal{H}_N(\kappa))}=\|\Op_{\hbar}^{w}(a)\|_{L^2(\mathbb{R}^{d})}.\end{equation}

\subsection{Quantization of toral automorphisms}

Let $A$ be a matrix in $Sp(2d,\mathbb{Z})$. As in the case of the hamiltonian flow on a manifold, we would like to quantize the dynamic associated to $A$ on the phase space, i.e. define a quantum progator associated to $A$. This can be done using the metaplectic `representation' of $Sp(2d,\mathbb{R})$~\cite{Fo} that defines for each matrix $A$ the unique (up to a phase) operator which satisfies:
$$\forall\rho\in\mathbb{R}^{2d},\ M(A)U_{\hbar}(\rho)M(A)^{-1}=U_{\hbar}(A\rho).$$
$M(A)$ is called the quantum propagator associated to $A$. It is a unitary operator on $L^2(\mathbb{R}^d)$ and it can be shown~\cite{BodB2}:
\begin{lemm}
Let $A$ be an element in $Sp(2d,\mathbb{Z})$ such that $1$ is not an eigenvalue of $A$. For each $N\in \mathbb{N}^{*}$, there exists at least one $\kappa_A\in[0,2\pi[^{2d}$ such that
$$M(A)\mathcal{H}_N(\kappa_A)=\mathcal{H}_N(\kappa_A).$$
$M_{\kappa_A}(A)$ denotes then the restriction of $M(A)$ to $\mathcal{H}_N(\kappa_A)$. It is an unitary operator.
\end{lemm}
Even if it was only stated for ergodic matrices, the proof of this lemma was given in~\cite{BodB2} (lemma $2.2$). The hypothesis that $1$ is not an eigenvalue is crucial in the proof of~\cite{BodB2}. We will say that an element $A$ in $Sp(2d,\mathbb{Z})$ is \emph{quantizable} if $1$ is not an eigenvalue of $A$. In the following, we will take $\kappa$ equal to the $\kappa_A$ (that also depends on $N$) given by this lemma. From all this, the following `exact' Egorov property can be shown for each $a\in\mathcal{C}^{\infty}(\mathbb{T}^{2d})$~\cite{Fo},~\cite{BdB}:
$$M_{\kappa}(A)^{-1}\Op_{\kappa}^{w}(a)M_{\kappa}(A)=\Op_{\kappa}^{w}(a\circ A).$$
\begin{rema} We underline that we do not need any assumption on $A$ (except that is is symplectic) to define $M(A)$ on $L^2(\mathbb{R}^d)$. In particular, for every $Q$ in $Sp(2d,\mathbb{R})$, we have that
$$\forall\ a\in S^0_{\nu}(1),\ M(Q)^{-1}\Op_{\hbar}^w(a)M(Q)=\Op_{\hbar}^w(a\circ Q).$$
\end{rema}

\subsection{Anti-Wick quantization}

The Weyl quantization has the nice properties that it satisfies an exact Egorov property and that for a symbol $a$, $\Op_{\kappa}^w(a)^*=\Op_{\kappa}^w(\overline{a})$. However, it does not satisfy the property that if $a$ is nonnegative then $\Op_{\kappa}^w(a)$ is also nonnegative. As our goal is to construct measures using a quantization procedure, we would like for simplicity to consider a positive quantization. This can be achieved by considering the anti-Wick quantization. To describe this quantization, we define the coherent state at point $0$ on $\mathbb{R}^d$:
$$|0\rangle(x):=\left(\frac{1}{\pi\hbar}\right)^{\frac{d}{4}}e^{-\frac{\|x\|^2}{2\hbar}}.$$
We define the translated coherent state at point $\rho\in\mathbb{R}^{2d}$ as $|\rho\rangle:=U(\rho)|0\rangle$. Using these coherent states, we can define a quantization procedure for a symbol $a$ in a nice class of symbol on $\mathbb{R}^{2d}$:
$$\Op_{\hbar}^{AW}(a):=\int_{\mathbb{R}^{2d}}a(\rho)|\rho\rangle\langle\rho|\frac{d\rho}{(2\pi\hbar)^d}.$$
It is obvious that this quantization is positive. It can be verified also that it satisfies the property of resolution of identity:
$$\Op_{\hbar}^{AW}(1)=\text{Id}_{L^2(\mathbb{R}^d)}=\int_{\mathbb{R}^{2d}}|\rho\rangle\langle\rho|\frac{d\rho}{(2\pi\hbar)^d}.$$
This quantization is related to the Weyl quantization. To see this, we can define the gaussian observable
$\tilde{G}_{\hbar}(x,\xi):=\frac{1}{(\pi\hbar)^{d}}e^{-\frac{\|x\|^2+\|\xi\|^2}{\hbar}}.$ For a bounded observable $a$, the relation between the two procedures of quantization is $\Op_{\hbar}^{AW}(a)=\Op_{\hbar}^{w}(a\star \tilde{G}_{\hbar}).$ Using Calder\'on-Vaillancourt theorem and the previous property, one can verify that $\|\Op_{\hbar}^{AW}(a)-\Op_{\hbar}^{w}(a)\|_{L^2(\mathbb{R}^d)}=O_a(\hbar).$ Now one can construct a positive quantization on the torus mimicking this positive quantization on $\mathbb{R}^{2d}$. To do this, we project the coherent states on the Hilbert space $\mathcal{H}_N(\kappa)$:
$$|\rho,\kappa\rangle:=S(\kappa)|\rho\rangle,$$
where $S(\kappa)$ is the projector defined by~(\ref{definition-projecteur-chat}). We define the anti-Wick quantization of an observable $a$ in $\mathcal{C}^{\infty}(\mathbb{T}^{2d})$ as follows:
$$\Op_{\kappa}^{AW}(a):=\int_{\mathbb{T}^{2d}}a(\rho)|\rho,\kappa\rangle\langle\rho,\kappa|\frac{d\rho}{(2\pi\hbar)^d}.$$
It satisfies that for a symbol $a$, $\Op_{\kappa}^{AW}(a)^*=\Op_{\kappa}^{AW}(\overline{a})$ and that the quantization is nonnegative. As in the case of the Weyl quantization, it is related to the quantization on $\mathbb{R}^{2d}$ by the integral representation~\cite{BdB}:
$$\Op_{\hbar}^{AW}(a)=\int_{[0,2\pi[^{2d}}\Op_{\kappa}^{AW}(a)d\kappa.$$
It also satisfies a resolution of identity property~\cite{BdB}:
$$\Op_{\kappa}^{AW}(1)=\text{Id}_{\mathcal{H}_N(\kappa)}=\int_{\mathbb{T}^{2d}}|\rho,\kappa\rangle\langle\rho,\kappa|\frac{d\rho}{(2\pi\hbar)^d}.$$
It is clear that $\|\Op_{\kappa}^{AW}(a)\|_{\mathcal{L}(\mathcal{H}_N(\kappa))}\leq\|a\|_{\infty}$.

\subsection{Semiclassical measures}

All these definitions allow to introduce the notion of semiclassical measures for the quantized cat-maps~\cite{BdB}:
\begin{def1}\label{d:semi-class-meas} Let $A$ be a matrix in $Sp(2d,\mathbb{Z})$ such that $1$ is not an eigenvalue of $A$. We call semiclassical measure of $(\mathbb{T}^{2d},A)$ any accumulation point of a sequence of measures of the form
$$\forall a\in\mathcal{C}^{\infty}(\mathbb{T}^{2d},\mathbb{C}),\ \tilde{\mu}^N(a):=\langle\psi^N|\Op_{\kappa}^{AW}(a)|\psi^N\rangle_{\mathcal{H}_N(\kappa)}=\int_{\mathbb{T}^{2d}}a(\rho)N\left|\left\langle\psi^N|\rho,\kappa\right\rangle_{\mathcal{H}_N(\kappa)}\right|^2d\rho,$$
where $(\psi^N)_N$ is a sequence of eigenvectors of $M_{\kappa}(A)$ in $\mathcal{H}_N(\kappa)$.
\end{def1}
The set of semiclassical measure defines a nonempty set of probability measures on the torus $\mathbb{T}^{2d}$. They are $A$-invariant measures using the following Egorov property:
\begin{prop}[Egorov property] Let $A$ be a matrix in $Sp(2d,\mathbb{Z})$ such that $1$ is not an eigenvalue of $A$. For every $a$ in $\mathcal{C}^{\infty}(\mathbb{T}^{2d})$, one has
$$\forall t \in\mathbb{R},\ M_{\kappa}(A)^{-t}\Op_{\kappa}^{AW}(a)M_{\kappa}(A)^{t}=\Op_{\kappa}^{AW}(a\circ A^t)+O_{a,t}(N^{-1}),
$$
where the constant involved in the remainder depends on $a$ and $t$.
\end{prop}
To conclude the presentation of our system, we underline that the set of semiclassical measures does not change if we consider another quantization procedure. For instance, we could have taken the Weyl procedure or any quantization $\Op_{\kappa}$ that satisfies, in the semiclassical limit,
$$\left\|\Op_{\kappa}(a)-\Op_{\kappa}^w(a)\right\|_{\mathcal{L}(\mathcal{H}_N(\kappa))}=O(N^{-\gamma}),$$
for some fixed positive $\gamma$.

\section{Symplectic linear algebra and Lyapunov exponents}\label{s:oseledets-sympl}

In this section, we collect some facts about symplectic matrices that we will use crucially in our proof. We refer the reader to chapter $1$ of~\cite{Lo} for more details. We fix a \emph{quantizable} matrix $A$ in $Sp(2d,\mathbb{Z})$, i.e. such that $1$ is not an eigenvalue of $A$. As theorem~\ref{t:main-theo} is trivial in the case where the spectrum is included in $\{z\in\mathbb{C}:|z|=1\}$, we also make the assumption that $A$ has an eigenvalue of modulus larger than $1$. We will denote
$$\lambda_{\max}=\sup\{\log|\beta|:\beta\ \text{is in the spectrum of}\ A\}.$$
\begin{rema} According to Kronecker's theorem ($2.5$ in~\cite{Na}), we know that if $A$ is an \emph{ergodic} matrix in $SL(2d,\mathbb{Z})$, then $\lambda_{\max}>0$.
\end{rema}
One can decompose $\mathbb{R}^{2d}$ into $A$-invariant subspaces called the stable, neutral and unstable spaces, i.e.
$$\mathbb{R}^{2d}:=E^-\oplus E^0\oplus E^+.$$
These subspaces satisfy various properties that we will use. The spectrum of the restriction of $A$ on the neutral space $E^0$ is included in $\{z\in\mathbb{C}:|z|=1\}$. The dimension of $E^0$ is even and we will denote it $2d_0$. The restriction of $A$ on the stable (resp. unstable) space $E^-$ (resp. $E^+$) has a spectrum included in $\{z\in\mathbb{C}:|z|<1\}$ (resp. $\{z\in\mathbb{C}:|z|>1\}$). These two subspaces have the same dimension equal to $d-d_0$ (which is by assumption positive). Moreover, there exist $r$ in $\mathbb{N}$ and $0<\lambda_1^+<\cdots<\lambda_r^+$ such that $E^+$ (resp. $E^-$) can be decomposed into $A$-invariant subspaces as follows:
$$E^+=E_1^+\oplus\cdots\oplus E_r^+\ \text{and}\ E^-=E_1^-\oplus\cdots\oplus E_r^-,$$
where the spectrum of the restriction of $A$ to $E^+_i$ (resp. $E^-_i$) is included in $\{z\in\mathbb{C}:|z|=e^{\lambda_i^+}\}$ (resp.  $\{z\in\mathbb{C}:|z|=e^{-\lambda_i^+}\}$). Moreover, one can verify that the subspaces $E^+_i$ and $E^-_i$ have the same dimension that we will denote $d_i$. The coefficients $\lambda_i^+$ are called the positive Lyapunov exponents of $A$. We underline that $\lambda_r^+=\lambda_{\max}$. With these notations, theorem~\ref{t:main-theo} can be rewritten that for any semiclassical measure $\mu$ associated to $A$, one has
\begin{equation}\label{e:main-theo}
h_{KS}(\mu,A)\geq\sum_{i=1}^rd_i\max\left(\lambda_i^+-\frac{\lambda_{\max}}{2},0\right).
\end{equation}
For the sake of simplicity, we will denote
$$\Lambda_+:=\sum_{i=1}^rd_i\lambda_i^+$$
and
$$\Lambda_0:=\sum_{i=1}^rd_i\max\left(\lambda_i^+-\frac{\lambda_{\max}}{2},0\right).$$
Our decomposition is exactly the Oseledets decomposition associated to the dynamical system $(\mathbb{T}^{2d},A,\mu)$. For our proof, we will need something stronger in order to apply tools of semiclassical analysis. Precisely, we will need a symplectic decomposition of $\mathbb{R}^{2d}$ into these subspaces. According to~\cite{Lo} (section $1.4$ to $1.7$), this decomposition is possible and we now recall the results from~\cite{Lo} that we will need. To do this, we introduce the $\diamond$-product of two matrices. Consider two real matrices $M_1$ in $M(2d',\mathbb{R})$ and $M_2$ in $M(2d'',\mathbb{R})$ of the block form
$$M_1:=\left(\begin{array}{cc}A_1&B_1\\C_1&D_1\end{array}\right)\ \text{and}\ M_2:=\left(\begin{array}{cc}A_2&B_2\\C_2&D_2\end{array}\right),$$
where $A_1$, $B_1$, $C_1$ and $D_1$ are in $M(d',\mathbb{R})$ and $A_2$, $B_2$, $C_2$ and $D_2$ are in $M(d'',\mathbb{R})$. The $\diamond$-product of $M_1$ and $M_2$ is defined as the following $2(d'+d'')$ matrix:
$$M_1\diamond M_2:=\left(\begin{array}{cccc}A_1&0&B_1&0\\0&A_2&0&B_2\\C_1&0&D_1&0\\0&C_2&0&D_2\end{array}\right).$$
We can use this product to rewrite our symplectic matrix $A$ in an adapted symplectic basis~\cite{Lo} (section $1.7$-theorem $3$). Precisely, for every $1\leq i\leq r$, one can construct an adapted $D_i$ in $Gl(d_i,\mathbb{R})$ such that the spectrum of $D_i$ is included in $\{z\in\mathbb{C}:|z|=e^{\lambda_i^+}\}$ and denote $A_i:=\text{diag}(D_i,D_i^{*-1})$ (setion $1.7$ in~\cite{Lo}). There exists also an adapted $A_0$ in $Sp(2d_0,\mathbb{R})$ such that the spectrum\footnote{We underline that $d_0$ will be equal to $0$ if all the Lyapunov exponents of $A$ are nonzero.} of $A_0$ is included in $\{z\in\mathbb{C}:|z|=1\}$ (section $1.5$ and $1.6$ in~\cite{Lo}). Using these matrices, it can be shown that there exists a symplectic matrix $Q$ in $Sp(2d,\mathbb{R})$ such that
\begin{equation}\label{e:diag-form}A=Q\left(A_0\diamond A_1\diamond\cdots A_r\right) Q^{-1}.\end{equation}
This tells us that we have a symplectic reduction adapted to the Oseledets decomposition. The results in~\cite{Lo} are more precise and we have only stated what we will need for our proof of theorem~\ref{t:main-theo}.

\section{Positive quantization adapted to the dynamic}\label{s:new-quant}

We have defined the set of semiclassical measures starting from the anti-Wick quantization. Yet, we have mentioned that this set does not depend on the choice of the quantization procedure. In this section, we will construct a new (positive) quantization procedure that is really adapted to the classical dynamic. To do this, we will mimick the construction of the anti-Wick quantization. In this case, we have seen that it corresponds to the Weyl quantization applied to the observable $a\star \tilde{G}_{\hbar}$. It means that we have made the convolution of the observable $a$ with a Gaussian observable which is localized in a ball of radius $\sqrt{\hbar}$.\\
Our strategy is to make a slightly different choice of function $G_{\hbar}$ which will be localized on an ellipsoid with lengths on each direction that depend on the Lyapunov exponent. For instance, if the Lyapunov exponent associated to the variable $(x_1,\xi_1)$ is larger than the one associated to the variable $(x_2,\xi_2)$, the ellipsoid will be larger in the second direction. We should also take care of not violating the uncertainty principle and as a consequence the radius of the ellipsoid will always be bounded from below by $\sqrt{\hbar}$.\\
In this section, we make this argument precise in the case of $\mathbb{R}^{2d}$ and then periodize the new quantization to get a quantization on the torus.

\subsection{An adapted convolution observable}

To construct our new new quantization on $\mathbb{R}^{2d}$, we introduce a Gaussian observable $G(x,\xi):=\exp\left(-\pi\|(x,\xi)\|^2\right),$ where $\|.\|$ is the euclidian norm on $\mathbb{R}^{2d}$. In the case of the anti-Wick quantization, we took the convolution of any bounded observable $a$ with $\displaystyle G\circ \left((\pi\hbar)^{-\frac{1}{2}}\text{Id}\right)$ to construct our quantization. Regarding the Oseledets decomposition of $A$ (see~(\ref{e:diag-form})), we would like to make a more strategical choice for the matrix we choose. To do this, we use the notations of section~\ref{s:oseledets-sympl} and for $\hbar>0$, we introduce a matrix $B(\hbar)$ of the following form:
$$B(\hbar):=Q\left(\begin{array}{cc}D_1(\hbar)&0\\0&D_1(\hbar)\end{array}\right)Q^{-1},$$
where $D_1(\hbar)$ is an element in $GL(d,\mathbb{R})$ of the form $$D_{1}(\hbar):=\left(\hbar^{-\frac{\epsilon_0}{2\lambda_{\max}}}\text{Id}_{\mathbb{R}^{d_0}},\hbar^{-\frac{\lambda_1^+}{2\lambda_{\max}}}\text{Id}_{\mathbb{R}^{d_1}},\cdots,\hbar^{-\frac{\lambda_r^+}{2\lambda_{\max}}}\text{Id}_{\mathbb{R}^{d_r}}\right).$$ In the previous definition, $\epsilon_0$ is some small fixed positive number that we keep fixed until the end of the proof. To simplify the expressions, we introduce the notation $\gamma_j^+:=\frac{\lambda_j^+}{2\lambda_{\max}}$ and $\gamma_0^+:=\frac{\epsilon_0}{2\lambda_{\max}}$. In particular, we have that $\|B(\hbar)^{-1}\|_{\infty}=\mathcal{O}(\hbar^{\gamma})$ for some fixed positive $\gamma$. Finally, we can define an `adapted' $\hbar$-Gaussian observable
$$G_{\hbar}:=2^{\frac{d}{2}}|\det B(\hbar)|^{\frac{1}{2}}G\circ B(\hbar).$$

\subsection{Positive quantization on $\mathbb{R}^{2d}$}\label{ss:positive-R}

We have constructed a convolution which seems to be adapted to the dynamic. We would also like to keep the nice `$u^*u$'-structure of the anti-Wick quantization. So, for a bounded observable $a$ in $\mathcal{C}^{\infty}(\mathbb{R}^{2d})$, we define
$$\Op_{\hbar}^+(a):=\Op_{\hbar}^w\left(a\star\left(G_{\hbar}\sharp G_{\hbar}\right)\right),$$
where $a\star b$ is the convolution product of two observables and $a\sharp b$ is the Moyal product of two observables (i.e. the symbol of $\Op_{\hbar}^w(a)\circ\Op_{\hbar}^w(b)$~\cite{DS}). We verify that
$$\Op_{\hbar}^+(a)=\int_{\mathbb{R}^{2d}}a(\rho_0)\Op_{\hbar}^w\left(\left(G_{\hbar}\sharp G_{\hbar}\right)(\bullet-\rho_0)\right)d\rho_0=\int_{\mathbb{R}^{2d}}a(\rho_0)\Op_{\hbar}^w\left(G_{\hbar}^{\rho_0}\right)^*\circ\Op_{\hbar}^w\left(G_{\hbar}^{\rho_0}\right)d\rho_0,$$
where $G_{\hbar}^{\rho_0}(\rho):=G_{\hbar}(\rho-\rho_0).$ So if $a\geq 0$, this defines a nonnegative operator. The following lemma says that $\Op_{\hbar}^+$ is a nice quantization procedure:
\begin{lemm}\label{l:equiv-procedure} Let $a$ be an observable in $S^{0}(1)$. We have
$$\left\|\Op_{\hbar}^w(a)-\Op_{\hbar}^+(a)\right\|_{L^2(\mathbb{R}^d)\rightarrow L^2(\mathbb{R}^d)}=O_a(\hbar^{\gamma}),$$
for some fixed positive $\gamma$ (depending only on $A$).
\end{lemm}
We postpone the proof of this lemma to appendix~\ref{a:eq-proc}. The strategy is the same as when one proves the equivalence of the anti-Wick quantization and the Weyl one. Precisely, we can prove that there exists an explicit kernel $K_{\hbar}(\rho_0)$ such that
$$a\star\left(G_{\hbar}\sharp G_{\hbar}\right)(\rho)=\int_{\mathbb{R}^{2d}}a(\rho+B(\hbar)^{-1}\rho_0)K_{\hbar}(\rho_0)d\rho_0,$$
where $\int_{\mathbb{R}^{2d}}K_{\hbar}(\rho_0)d\rho_0=1$.

\subsection{Periodization of observables}

We have just defined a new quantization procedure on $\mathbb{R}^{2d}$ which is related to the classical dynamic associated to the matrix $A$. To study our problem, we need to restrict this quantization procedure to $\mathcal{H}_{N}(\kappa)$. To do this, we define
\begin{equation}\Op_{\kappa}^+(a):=\Op_{\kappa}^w\left(a\star\left(G_{\hbar}\sharp G_{\hbar}\right)\right).\end{equation}
Thanks to lemma~\ref{l:equiv-procedure} and to the decomposition of $L^2(\mathbb{R}^{d})$ along the spaces $\mathcal{H}_N(\kappa)$, we know that $\|\Op_{\kappa}^+(a)-\Op_{\kappa}^w(a)\|_{\mathcal{L}(\mathcal{H}_N(\kappa))}=O_a(\hbar^{\gamma})$ (see also~\cite{ReSi} (theorem XIII.83)). The explicit form of this procedure is given by
$$\Op_{\kappa}^+(a)=\sum_{r\in\mathbb{Z}^{2d}}\left(\int_{\mathbb{T}^{2d}}e^{2\imath\pi\langle \rho,Jr\rangle}\int_{\mathbb{R}^{2d}}a(\rho_0)\left(G_{\hbar}\sharp G_{\hbar}\right)(\rho-\rho_0)d\rho_0d\rho\right)U_{\hbar}\left(\frac{r}{N}\right).
$$
For our purpose, we would like to verify that it remains a positive quantization procedure with a nice structure. To see this, we introduce the following periodization operators on $\mathcal{S}(\mathbb{R}^{2d})$:
\begin{equation}\label{e:trans-op}\forall\ \rho\in\mathbb{R}^{2d},\ \forall\ F\in\mathcal{S}(\mathbb{R}^{2d}),\ T_{\rho}(F)(\rho'):=\sum_{r\in\mathbb{Z}^{2d}}F\left(\rho'+r-\frac{J\rho}{2N}\right)e^{2\imath\pi\langle r+\rho',\rho\rangle}.\end{equation}
In particular, $T_0(G_{\hbar}^{\rho_0}\sharp G_{\hbar}^{\rho_0})$ is an element in $\mathcal{C}^{\infty}(\mathbb{T}^{2d})$ and it allows us to rewrite
$$\Op_{\kappa}^+(a)=\int_{\mathbb{T}^{2d}}a(\rho_0)\Op_{\kappa}^w\left(T_0(G_{\hbar}^{\rho_0}\sharp G_{\hbar}^{\rho_0})\right)d\rho_0.$$
The translation operators $T_{\rho}$ satisfy the following property:
\begin{prop}\label{p:new-quant-torus} Let $F_1$ and $F_2$ be two elements in $\mathcal{S}(\mathbb{R}^{2d})$. One has
$$\Op_{\hbar}^w(T_0(\overline{F}_1\sharp F_2))=\int_{\mathbb{T}^{2d}}\Op_{\hbar}^w(T_{\rho}F_1)^*\circ\Op_{\hbar}^{w}(T_{\rho}F_2)d\rho.$$
\end{prop}
We postpone the proof of this lemma (which is just a careful application of the Poisson formula) to appendix~\ref{s:technical-lemmas}. This proposition provides an alternative form for our quantization procedure, i.e.
\begin{equation}\forall\kappa\in[0,2\pi[^{2d},\ \Op_{\kappa}^+(a)=\int_{\mathbb{T}^{2d}}a(\rho_0)\int_{\mathbb{T}^{2d}}\Op_{\kappa}^w(T_{\rho}(G_{\hbar}^{\rho_0}))^*\circ\Op_{\kappa}^w(T_{\rho}(G_{\hbar}^{\rho_0}))d\rho d\rho_0.\end{equation}
In particular, it implies that $\Op_{\kappa}^+$ is a nonnegative quantization procedure. We also underline that we have the following resolution of identity:
\begin{equation}\label{resolution-id}
\text{Id}_{\mathcal{H}_N(\kappa)}=\int_{\mathbb{T}^{2d}}\int_{\mathbb{T}^{2d}}\Op_{\kappa}^w(T_{\rho}(G_{\hbar}^{\rho_0}))^*\circ\Op_{\kappa}^w(T_{\rho}(G_{\hbar}^{\rho_0}))d\rho d\rho_0.
\end{equation}
\begin{rema} These last two formulas are the analogues of the ones obtained for the anti-Wick quantization. The expressions seems more complicated but we will see that it is more adapted to the dynamic induced by $A$.
\end{rema}

\subsection{Long times Egorov property} In this last paragraph, we show that as the anti-Wick procedure, the quantization procedure $\Op_{\kappa}^+$ satisfies an Egorov property until times of order $T_E(N):=\frac{\log N}{2\lambda_{\max}}$. We fix some positive $\epsilon\ll\min(\epsilon_0,\lambda_1)$ and define the Ehrenfest time
\begin{equation}m_E(N):=\left[\frac{1-\epsilon}{2\lambda_{\max}}\log N\right],\end{equation}
The parameter $\epsilon$ will be kept fixed (until the end of the proof of theorem~\ref{t:main-theo}). In order to state our result, we denote $\mu^N$ the measure associated to the unit eigenvector $\psi_N$, i.e.
$$\mu^N(a):=\langle\psi_N|\Op_{\kappa}^{+}(a)|\psi_N\rangle_{\mathcal{H}_N(\kappa)}=\int_{\mathbb{T}^{2d}}a(\rho_0)\int_{\mathbb{T}^{2d}}\left\|\Op_{\kappa}^w(T_{\rho}(G_{\hbar}^{\rho_0}))\psi_N\right\|_{\mathcal{H}_N(\kappa)}^2d\rho d\rho_0.$$
One can show the following (pseudo)-invariance property of the measures $\mu^N$ until time $m_E(N)$:
\begin{prop}\label{non-exact-egorov-chat} Let $(\psi^N)_N$ be a sequence of unit eigenvectors of $M_{\kappa}(A)$ in $\mathcal{H}_N(\kappa)$ and $\mu^N$ the associated sequence of measures. Then, for every positive $\epsilon$, one has
\begin{equation}\forall a\in\mathcal{C}^{\infty}(\mathbb{T}^2,\mathbb{C}),\ \forall |t|\leq m_E(N),\ \mu^N(a\circ A^t)=\mu^N(a)+o_{a,\epsilon}(1),\end{equation}
where the constant in remainder depends only on $a$ and $\epsilon$.
\end{prop}
\emph{Proof.} We have an exact Egorov property for the Weyl quantization. In particular, it tells us that, for every integer $t$,
$$\left\|\Op_{\kappa}^+(a\circ A^t)-\Op_{\kappa}^+(a)(t)\right\|_{\mathcal{L}(\mathcal{H}_N(\kappa))}=\left\|\Op_{\kappa}^+(a\circ A^t)-\Op_{\kappa}^w(a\circ A^t)\right\|_{\mathcal{L}(\mathcal{H}_N(\kappa))}+O_a(\hbar^{\gamma}),$$
where $\Op_{\kappa}^+(a)(t):=M_{\kappa}(A)^{-t}\Op_{\kappa}^+(a)M_{\kappa}(A)^{t}$. From the decomposition of the space $L^2(\mathbb{R}^d)$ along the spaces $\mathcal{H}_N(\kappa)$, we know that
$$\left\|\Op_{\kappa}^+(a\circ A^t)-\Op_{\kappa}^w(a\circ A^t)\right\|_{\mathcal{L}(\mathcal{H}_N(\kappa))}\leq\left\|\Op_{\hbar}^+(a\circ A^t)-\Op_{\hbar}^w(a\circ A^t)\right\|_{\mathcal{L}(L^2(\mathbb{R}^d))}.$$
Recall that we know that, for a bounded symbol $b$, $\Op_{\hbar}^+(b)$ is equal to the operator $\Op_{\hbar}^w(b\star\left(G_{\hbar}\sharp G_{\hbar}\right))$ and that $b\star\left(G_{\hbar}\sharp G_{\hbar}\right)(\rho)=\int_{\mathbb{R}^{2d}}b(\rho+B(\hbar)^{-1}\rho_0)K_{\hbar}(\rho_0)d\rho_0$ (see paragraph~\ref{ss:positive-R} and appendix~\ref{a:eq-proc}). We write this formula for $b=a\circ A^t$ and combine it with the Taylor formula. We find that
$$(a\circ A^t)\star\left(G_{\hbar}\sharp G_{\hbar}\right)(\rho)=a\circ A^t(\rho)+\int_{\mathbb{R}^{2d}}K_{\hbar}(\rho_0)\int_0^1\left(d_{\rho+sB(\hbar)^{-1}\rho_0}a\right).(A^tB(\hbar)^{-1})\rho_0dsd\rho_0.$$
Appendix~\ref{a:eq-proc} gives us an exact expression for $K_{\hbar}$. We can compute the derivatives of the second term of the sum and according to Calder\'on-Vaillancourt theorem (see paragraph~\ref{ss:weyl-quant}), we finally find that
\begin{equation}\label{e:new-Egorov-property}\left\|\Op_{\kappa}^+(a\circ A^t)-M_{\kappa}(A)^{-t}\Op_{\kappa}^+(a)M_{\kappa}(A)^{t}\right\|_{\mathcal{L}(\mathcal{H}_N(\kappa))}=O_a\left(\|A^tB(\hbar)^{-1}\|_{\infty}\right).
\end{equation}
As $B(\hbar)$ was constructed to be adapted to the classical dynamic induced by $A$, we can verify that this last equality allows to conclude the proof of proposition~\ref{non-exact-egorov-chat}.$\square$

\section{Proof of theorem~\ref{t:main-theo}}\label{s:proof-maintheo}

We consider a semiclassical measure $\mu$. Without loss of generality, we can suppose that it is constructed from $\Op_{\kappa}^+$ and that it is associated to a sequence of eigenvectors $\psi_{N_k}$ of $M_{\kappa}(A)$ in $\mathcal{H}_{N_k}(\kappa)$ where $(N_{k})_k$ is an increasing sequence of integers. Precisely, we have $$\forall\ a\in\mathcal{C}^{\infty}(\mathbb{T}^2,\mathbb{C}),\ \mu(a)=\lim_{k\rightarrow+\infty}\langle\psi_{N_k}|\Op_{\kappa}^{+}(a)|\psi_{N_k}\rangle_{\mathcal{H}_{N_k}(\kappa)}.$$
We recall that we have denoted $\mu^{N_k}(a)=\langle\psi_{N_k}|\Op_{\kappa}^{+}(a)|\psi_{N_k}\rangle_{\mathcal{H}_{N_k}(\kappa)}$. To simplify notations, we will not mention $k$ in the following of this article. We start our proof by fixing a finite measurable partition $\mathcal{Q}$ of small diameter $\delta_0$ whose boundary is not charged\footnote{The parameter $\delta_0$ is small and fixed for all the article: it has no vocation to tend to $0$.} by $\mu$ (paragraph $2.2.8$ in~\cite{AN2}). We denote $\eta(x):=-x\log x$ (with the convention $0\log 0=0$). We recall that the Kolmogorov-Sinai entropy of the measure $\mu$ for the partition $\mathcal{Q}$ can be defined as~\cite{Wa}
$$h_{KS}(\mu,A,\mathcal{Q}):=\lim_{m\rightarrow+\infty}\frac{1}{2m}\sum_{|\alpha|=2m}\eta\left(\mu\left(A^mQ_{\alpha_{-m}}\cdots\cap A^{-(m-1)}Q_{\alpha_{m-1}}\right)\right),$$
where $\alpha_j$ varies in $\{1,\cdots,K\}$ ($K$ is the cardinal of $\mathcal{Q}$).

\subsection{Using the entropic uncertainty principle}

Our quantization is defined for smooth observables on the torus. So we start by defining a smooth partition $(P_i)_{i=1}^K$ of observables in $\mathcal{C}^{\infty}(\mathbb{T}^2,[0,1])$ (of small support of diameter less than $2\delta$) that satisfies the following property of partition of $\mathbb{T}^{2d}$:
\begin{equation}\label{partition-identité-tore-classique}
\forall\rho\in\mathbb{T}^{2d},\ \sum_{i=1}^KP_i^2(\rho)=1.
\end{equation}
Mimicking the definition of Kolmogorov-Sinai entropy, we define the quantum entropy of $\psi_N$ with respect to $\mathcal{P}$:
\begin{equation}\label{définition-entropie-quantique-tore}h_{2m}(\psi_N,\mathcal{P}):=-\sum_{|\alpha|=2m}\mu^N(\mathbf{P}_{\alpha}^2)\log\mu^N(\mathbf{P}_{\alpha}^2),
\end{equation}
where $\mathbf{P}_{\alpha}:=\prod_{j=-m}^{m-1}P_{\alpha_j}\circ A^j$ for $\alpha:=(\alpha_{-m},\cdots,\alpha_{m-1})$. One can verify that for fixed $m$, we have
\begin{equation}\label{passage-limite-semi-classique-entropie}
h_{2m}(\mu,\mathcal{P}):=-\sum_{|\alpha|=2m}\mu(\mathbf{P}_{\alpha}^2)\log\mu(\mathbf{P}_{\alpha}^2)=\lim_{N\rightarrow\infty}h_{2m}(\psi_N,\mathcal{P}).
\end{equation}
So, for a fixed $m$ and as $N\rightarrow\infty$, the quantum entropy we have just defined tends to the usual entropy of $\mu$ at time $2m$ (with the notable difference that we consider smooth partitions). Our crucial observation to apply the entropic uncertainty principle is that we have the following partition of identity for $\mathcal{H}_N(\kappa)$:
\begin{equation}\label{partition-fondamentale}
\sum_{|\alpha|=2m}\int_{\mathbb{T}^{2d}}\int_{\mathbb{T}^{2d}}\mathbf{P}_{\alpha}^2(\rho_0)\Op_{\kappa}^w(T_{\rho}(G_{\hbar}^{\rho_0}))^*\circ\Op_{\kappa}^w(T_{\rho}(G_{\hbar}^{\rho_0}))d\rho d\rho_0=\text{Id}_{\mathcal{H}_N(\kappa)}.
\end{equation}
This partition of identity is derived from equation~(\ref{resolution-id}) and is crucial to apply the entropic uncertainty principle~\ref{ppe-incertitude-partition-chat}. In fact, it can be applied for $\mathcal{H}=L^2(\mathbb{T}^{4d},\mathcal{H}_N(\kappa))$ and $\tilde{\mathcal{H}}=\mathcal{H}_{N}(\kappa)$. For $\rho$ in $\mathbb{T}^{2d}$ and $\psi$ in $\mathcal{H}_{N}(\kappa)$, we define $$\mathbf{\pi}_{\alpha}|\psi\rangle(\rho,\rho_0):=P_{\alpha}(\rho_0)\Op_{\kappa}^w(T_{\rho}(G_{\hbar}^{\rho_0}))|\psi\rangle.$$
This defines a linear application from $\mathcal{H}_{N}(\kappa)$ to $L^2(\mathbb{T}^{4d},\mathcal{H}_N(\kappa))$ and its adjoint is given by
$$\pi_{\alpha}^{\dagger}f:=\int_{\mathbb{T}^{4d}}P_{\alpha}(\rho_0)\Op_{\kappa}^w(T_{\rho}(G_{\hbar}^{\rho_0}))^*f(\rho,\rho_0) d\rho d\rho_0,$$ for $f$ in $L^2(\mathbb{T}^{4d},\mathcal{H}_N(\kappa))$. It defines a quantum partition of identity as it satisfies the relation
$\sum_{|\alpha|=2m}\pi_{\alpha}^{\dagger}\pi_{\alpha}=\text{Id}_{\mathcal{H}_N(\kappa)}.$ Applying the entropic uncertainty principle for this partition and $U=M_{\kappa}(A)^n$, we bound $\|\pi_{\alpha}M_{\kappa}(A)^n\pi_{\beta}^{\dagger}\|_{\mathcal{L}(L^2(\mathbb{T}^{4d},\mathcal{H}_N(\kappa)))}$ and derive the following corollary:
\begin{coro}\label{c:entropic} Using the previous notations, one has
$$\forall n\in\mathbb{N},\ \forall m\in\mathbb{N},\ h_{2m}(\psi_N,\mathcal{P})\geq-\log\sup_{|\alpha|=2m}\{\text{Leb}(\mathbf{P}_{\alpha}^2)\}-\log c(A,n),$$
where
$\displaystyle c(A,n):=\sup_{\rho,\rho',\rho_0,\rho_0'\in\mathbb{T}^{2d}}\left\{\left\|\Op_{\kappa}^w(T_{\rho}(G_{\hbar}^{\rho_0}))M_{\kappa}(A)^n\Op_{\kappa}^w(T_{\rho'}(G_{\hbar}^{\rho_0'}))^*\right\|_{\mathcal{L}(\mathcal{H}_N(\kappa))}\right\}.$
\end{coro}
\begin{rema} If we had considered the anti-Wick quantization, we would have used the fundamental relation
$$\sum_{|\alpha|=2m}\int_{\mathbb{T}^{2d}}\mathbf{P}_{\alpha}^2(\rho_0)|\rho_0,\kappa\rangle\langle\kappa,\rho_0|d\rho_0=\text{Id}_{\mathcal{H}_N(\kappa)}.$$
In this case~\cite{GR3}, we would obtain a similar corollary as~\ref{c:entropic}. In particular, the quantity $c(A,n)$ in the lower bound would be replaced by the $\sup_{\rho_0,\rho_0'}N^d\left|\langle\kappa,\rho'_0|M_{\kappa}(A)^{t}|\rho_0,\kappa\rangle_{\mathcal{H}_N(\kappa)}\right|$. We will see that it is not sufficient to deduce our main theorem.
\end{rema}

\subsection{Estimate of $c(A,n)$}
In section~\ref{main-estimate}, we will prove the following theorem:
\begin{theo}\label{t:main-est} Let $A$ be a \emph{quantizable} matrix and let $\epsilon$ be some (small) positive number. For every positive $\delta$ (small enough), there exists a constant $C$ such that, for $n:=n_E(\hbar)=[(1-\epsilon)|\log\hbar|/\lambda_{\max}]$,
$$c(A,n)\leq C|\det B(\hbar)|\hbar^{-\delta-\frac{\Lambda_+}{\lambda_{\max}}\epsilon} e^{-n_{E}(\hbar)\Lambda_0},$$
where $\Lambda_+$ and $\Lambda_0$ depend on the Lyapunov exponents of $A$ and were defined in section~\ref{s:oseledets-sympl}.
\end{theo}
Recall that we have the relation $2\pi\hbar N=1$. If we consider $\delta\ll\epsilon$, then the quantum entropy at time $2m$ is bounded from below as follows:
\begin{equation}\label{inégalité-temps-Ehrenfest-entropie}
\forall m\geq 1,\ h_{2m}(\psi_N,\mathcal{P})\geq\frac{\log N}{\lambda_{\max}}\left((\Lambda_0-\Lambda_+)(1-2\epsilon)-d_0\epsilon_0\right)-\log\sup_{|\alpha|=2m}\{\text{Leb}(\mathbf{P}_{\alpha}^2)\}+\tilde{C}.
\end{equation}
This estimate is our main simplification compared with~\cite{AN2} as it will only use estimates of gaussian integrals. We underline that this theorem plays a crucial role in our proof\footnote{The other term on the lower bound will be estimate thanks to the computation of the entropy of the Lebesgue measure.} as it replaces all the discussion of section $3$ in~\cite{AN2}.
\begin{rema} Even if we postpone the proof of theorem~\ref{t:main-est} to the next section, we can explain how it is related to previous results on the quantum correlation function~\cite{BodB},~\cite{FNdB}. In fact, it we had used the anti-Wick quantization in our application of the entropic uncertainty principle, we would have had to give an estimate on the quantum correlation function
$$N^d\left|\langle\kappa,\rho'|M_{\kappa}(A)^{t}|\rho,\kappa\rangle_{\mathcal{H}_N(\kappa)}\right|\leq N^d\sum_{r\in\mathbb{Z}^{2d}}\left|\left\langle 0\left|M_{\kappa}(A)^{-\frac{t}{2}}U_{\hbar}\left(r+A^{-\frac{t}{2}}\rho'+A^{\frac{t}{2}}\rho\right)M_{\kappa}(A)^{\frac{t}{2}}\right|0\right\rangle_{L^2(\mathbb{R}^d)}\right|.$$
In~\cite{BodB},~\cite{FNdB}, it was proved that as long as $|t|\leq\frac{1-\epsilon}{\lambda_{\max}}\log N$, the only term that contributes to the previous sum is the term $(0,0)$, precisely:
$$\forall|t|\leq 2m_E(N),\ \forall\rho\in\mathbb{T}^{2d},\ \forall\rho'\in\mathbb{T}^{2d},\ N^d\left|\langle\kappa,\rho'|M_{\kappa}(A)^{t}|\rho,\kappa\rangle_{\mathcal{H}_N(\kappa)}\right|\leq CN^d e^{-\frac{\Lambda_+t}{2}}.$$
Using this inequality, we would only find $\frac{3}{2}\Lambda_+-d\lambda_{\max}$ in the lower bound of theorem~\ref{t:main-theo} (as it was obtained in~\cite{AN2}). This explains why we have introduced this new quantization. In fact, with our new choice, we are able to obtain a similar bound (with a similar method) but the correction is not anymore $N^d$.
\end{rema}

\subsection{Subadditivity of the quantum entropy}
Now, we have to find a time $m$ for which inequality~(\ref{inégalité-temps-Ehrenfest-entropie}) is optimal. It will depend on $N$ and the last difficulty is that if $m(N)$ grows too fast with $N$, $h_{2m(N)}(\psi_N,\mathcal{P})$ has no particular reason to tend to $h_{KS}(\mu,A)$ in the semiclassical limit. We have to be careful and we first verify that classical arguments from ergodic theory (subadditivity of the entropy) can be adapted for the quantum entropy as long as $m\leq \log N/(2\lambda_{\max})$. In particular, we prove that the sequence $\frac{1}{2m_0}h_{2m_0}(\psi_N,\mathcal{P})$ is `almost' decreasing until the Ehrenfest time (see appendix~\ref{proof-lemma}):
\begin{lemm}\label{lemme-chat} We denote $m_E(N)=[(1-\epsilon)\log N/(2\lambda_{\max})]$ and fix an integer $m_0$. We have then
$$\frac{1}{2m_E(N)}h_{2m_E(N)}(\psi_N,\mathcal{P})\leq \frac{1}{2m_0}h_{2m_0}(\psi_N,\mathcal{P})+R(m_0,N),$$
where $R(m_0,N)$ is a remainder that satisfies $\forall m_0\in\mathbb{N},\ \lim_{N\rightarrow\infty}R(m_0,N)=0.$
\end{lemm}
Combining this lemma with the entropic estimation~(\ref{inégalité-temps-Ehrenfest-entropie}), we have, for every fixed $m_0>0$,
\begin{equation}\label{temps-m0-chat} \frac{1}{2m_0}h_{2m_0}(\psi_N,\mathcal{P})+\tilde{R}(m_0,N)\geq\left((\Lambda_0-\Lambda_+)-\frac{d_0\epsilon_0}{1-2\epsilon}\right)-\frac{1}{2m_E(N)}\log\sup_{|\alpha|=2m_E(N)}\{\text{Leb}(\mathbf{P}_{\alpha}^2)\}.\end{equation}
where $\tilde{R}(m_0,N)$ is a remainder that satisfies $\forall m_0\in\mathbb{N},\ \lim_{N\rightarrow\infty}\tilde{R}(m_0,N)=0.$

\subsection{The conclusion}
To conclude, it remains to bound the quantity $\sup_{|\alpha|=2m_E(N)}\{\text{Leb}(\mathbf{P}_{\alpha}^2)\}$. To do this, we underline that, for each $\alpha$ of length $2m$,
$$\forall x\in \text{supp}(\mathbf{P}_{\alpha}^2),\ \text{Leb}(\mathbf{P}_{\alpha}^2)\leq\text{Leb}(\text{supp}(\mathbf{P}_{\alpha}^2))\leq\text{Leb}(B(x,2\delta_0,2m)),$$
where $B(x,2\delta_0,2m):=\{y\in\mathbb{T}^2:\forall j\in[-m,m-1],\ d(A^jx,A^jy)<2\delta_0\}$, where $d$ is the metric induced on $\mathbb{T}^{2d}$ by the Euclidean norm on $\mathbb{R}^{2d}$. By induction and using the invariance of the metric $d$, we know that for every $x$ in $\mathbb{T}^{2d}$ and for every $k$ in $\mathbb{Z}$, $A^{-k}B(A^kx,2\delta_0)=x+A^{-k}B(0,2\delta_0)$. Then, using the invariance by translation of the Lebesgue measure, we know that for every $x$ in $\mathbb{T}^{2d}$, $\text{Leb}(B(x,2\delta_0,2m))=\text{Leb}(B(0,2\delta_0,2m))$. Combining~\cite{BKa} and theorem $8.15$ from~\cite{Wa}, we know that $\text{Leb}(B(0,2\delta_0,2m))\leq C_{\delta_0}e^{-2m(\Lambda_+-\epsilon)}$. We use this last inequality and we make $N$ tends to infinity in~(\ref{temps-m0-chat}). It gives, for every positive $m_0$,
$$\frac{1}{2m_0}h_{2m_0}(\mu,\mathcal{P})\geq\Lambda_0-\frac{d_0\epsilon_0}{1-2\epsilon}-\epsilon.$$
This last inequality holds for all (small enough) smoothing $\mathcal{P}$ of the partition $\mathcal{Q}$. The lower bound does not depend on the derivatives of $\mathcal{P}$ so we we can replace the smooth partition $\mathcal{P}$ by the true partition $\mathcal{Q}$ in the definition of $h_{2m_0}(\mu,\mathcal{P})$. We let $m_0$ tends to infinity, then $\epsilon$ to $0$ and finally $\epsilon_0$ to $0$. We find
$$h_{KS}(\mu,A)\geq h_{KS}(\mu,A,\mathcal{Q})\geq\Lambda_0.\square$$

\section{The main estimate: proof of theorem~\ref{t:main-est}}\label{main-estimate}

In this section, we want to prove theorem~\ref{t:main-est}, i.e. give an estimate of $c(A,n)$. We underline that the spirit of the proof will  be similar to the proof of estimates on the propagation of coherent states under the quantum propagator~\cite{BodB},~\cite{FNdB}.\\
First, we use exact Egorov property for the Weyl quantization and we find
$$c(A,n):=\sup_{\rho,\rho',\rho_0,\rho_0'\in\mathbb{T}^{2d}}\left\{\left\|\Op_{\kappa}^w\left(T_{\rho}(G_{\hbar}^{\rho_0})\circ A^{\frac{n}{2}}\right)\Op_{\kappa}^w\left(\overline{T_{\rho'}(G_{\hbar}^{\rho_0'})}\circ A^{-\frac{n}{2}}\right)\right\|_{\mathcal{L}(\mathcal{H}_N(\kappa))}\right\}.$$
As $\Op_{\kappa}^w(a)$ is the restriction of $\Op_{\hbar}^w(a)$ to $\mathcal{H}_N(\kappa)$ and using the decomposition of $L^2(\mathbb{R}^d)$ along the $\mathcal{H}_N(\kappa)$, we know that
\begin{equation}\label{e:c(A,n)}c(A,n)\leq\sup_{\rho,\rho',\rho_0,\rho_0'\in\mathbb{T}^{2d}}\left\{\left\|\Op_{\hbar}^w\left(T_{\rho}\left(G_{\hbar}^{\rho_0}\right)\circ A^{\frac{n}{2}}\right)\Op_{\hbar}^w\left(\overline{T_{\rho'}(G_{\hbar}^{\rho_0'})}\circ A^{-\frac{n}{2}}\right)\right\|_{\mathcal{L}(L^2(\mathbb{R}^d))}\right\}.\end{equation}
Our goal is to give an estimate on this last quantity. Let us explain what our strategy will be. We underline that the symbols of the two pseudodifferentials operators are defined as infinite sums over $\mathbb{Z}^{2d}$. We will proceed in three steps. First, we will prove that the product of the two terms centered in $(0,0)$ satisfy the bound we expect (section~\ref{s:leadterm}). Then, we will prove that products of two terms not centered at the same element of $\mathbb{Z}^{2d}$ are negligible (section~\ref{s:negligible}). Finally, we will combine these two estimates with the Coltar Stein theorem~\cite{DS} in order to conclude (section~\ref{s:Coltar}).

\subsection{The leading term}\label{s:leadterm} We start our estimate by giving a bound on the term centered on $(0,0)$ in $\mathbb{Z}^{4d}$. It means that we will look at the norm of the operator
\begin{equation}\label{e:lead-op}\Op_{\hbar}^w\left(G_{\hbar}\left( A^{\frac{n}{2}}\bullet-\pi\hbar J\rho-\rho_0\right)e^{2\imath\pi\langle A^{\frac{n}{2}}\bullet-\rho_0|\rho\rangle}\right)\Op_{\hbar}^w\left(G_{\hbar}\left( A^{-\frac{n}{2}}\bullet-\pi\hbar J\rho'-\rho_0'\right)e^{2\imath\pi\langle A^{-\frac{n}{2}}\bullet-\rho_0'|\rho'\rangle}\right)^*,\end{equation}
where $\rho,\rho',\rho_0,\rho_0'\in\mathbb{T}^{2d}$. Precisely, we will prove in this paragraph the following proposition:
\begin{prop}\label{p:leadterm} Let $\rho$, $\rho'$, $\rho_0$ and $\rho_0'$ be elements in $\mathbb{T}^{2d}$. Let $\epsilon$ be some (small) fixed positive number. Then, for every positive $\delta$ and for $n=n_E(\hbar):=[(1-\epsilon)|\log\hbar|/\lambda_{\max}]$, one has
$$\|(\ref{e:lead-op})\|_{\mathcal{L}(L^2(\mathbb{R}^d))}\leq C|\det B(\hbar)|\hbar^{-\delta-\epsilon\frac{\Lambda_+}{\lambda_{\max}}}\exp\left(-\sum_{i:2\lambda_i^+-\lambda_{\max}>0}d_i\left(\lambda_i^+-\frac{\lambda_{\max}}{2}\right)n_E(\hbar)\right),$$
where $\Lambda_+:=\sum_{i=0}^rd_i\lambda_i^+$ and where the constant $C$ is uniform for $\rho$, $\rho'$, $\rho_0$ and $\rho_0'$ in $\mathbb{T}^{2d}$
\end{prop}
\begin{rema} We underline that this estimate is exactly the one of theorem~\ref{t:main-est}. We will check in the following paragraphs that terms not centered in $(0,0)$ have a smaller contribution. We also underline that we restrict ourselves to the time $n_E(\hbar)$. In the following sections (\ref{s:negligible} and~\ref{s:Coltar}), we will verify that the main contribution comes from the term centered in $(0,0)$ only if we restrict ourselves to $0\leq n\leq n_E(\hbar)$. Moreover, it will be clear in the proof that the bound is the smallest possible for $n=n_E(\hbar)$ if we only consider the range of times $0\leq n\leq n_E(\hbar)$.
\end{rema}

\subsubsection{First observations}
For simplicity of notations, we introduce two auxiliary matrices
$$A_+(n,\hbar):=Q^{-1}\sqrt{\hbar}B(\hbar)A^{\frac{n}{2}}Q\ \text{and}\ A_-(n,\hbar):=Q^{-1}\sqrt{\hbar}B(\hbar)A^{-\frac{n}{2}}Q.$$
Recall that, using the notations of section~\ref{s:oseledets-sympl}, we have
$$A_+(n,\hbar)=\left(\hbar^{\frac{\lambda_{\max}-\epsilon_0}{2\lambda_{\max}}}A_0^{\frac{n}{2}}\right)\diamond\left(\hbar^{\frac{\lambda_{\max}-\lambda_1^+}{2\lambda_{\max}}}A_1^{\frac{n}{2}}\right)\cdots\diamond\left(\hbar^{\frac{\lambda_{\max}-\lambda_r^+}{2\lambda_{\max}}}A_r^{\frac{n}{2}}\right)$$
and
$$A_-(n,\hbar)=\left(\hbar^{\frac{\lambda_{\max}-\epsilon_0}{2\lambda_{\max}}}A_0^{-\frac{n}{2}}\right)\diamond\left(\hbar^{\frac{\lambda_{\max}-\lambda_1^+}{2\lambda_{\max}}}A_1^{-\frac{n}{2}}\right)\cdots\diamond\left(\hbar^{\frac{\lambda_{\max}-\lambda_r^+}{2\lambda_{\max}}}A_r^{-\frac{n}{2}}\right).$$
We would like now to write the norm we have to estimate into a simpler form using these new notations. First, we underline that if we define
$$V_{\hbar}u(x):=\hbar^{\frac{d}{4}}u(\sqrt{\hbar} x),$$
then, for any bounded operators $\Op_{\hbar}^w(a)$ and $\Op_{\hbar}^w(b)$ on $L^2(\mathbb{R}^d)$, one has
$$\left\|\Op_{\hbar}^w(a)\Op_{\hbar}^w(b)\right\|_{\mathcal{L}(L^2(\mathbb{R}^d))}=\left\|\Op_{1}^w(a\circ(\sqrt{\hbar}\text{Id}_{\mathbb{R}^{2d}}))\Op_{1}^w(b\circ(\sqrt{\hbar}\text{Id}_{\mathbb{R}^{2d}}))\right\|_{\mathcal{L}(L^2(\mathbb{R}^d))}.$$
Moreover, we can use that the matrix $Q$ is an element of $Sp(2d,\mathbb{R})$. In particular, its quantization $M(Q)$ (via the metaplectic representation) satisfies an exact Egorov property~\cite{DS}. Using these two observations and defining $\tilde{\Gamma}_{\rho}(w):=G\circ Q(w)e^{2\imath\pi\langle w|\rho\rangle}$, we can deduce that the norm of the operator~(\ref{e:lead-op}) is bounded by
\begin{equation}
2^d|\det B(\hbar)|\sup_{\rho,\rho'\in[-M,M]^{2d};\rho_0,\rho_0'\in\mathbb{R}^{2d}}\left\|\Op_{1}^w\left(\tilde{\Gamma}_{\rho}\left(A_+(n,\hbar)\bullet-\rho_0\right)\right)\Op_{1}^w\left(\tilde{\Gamma}_{\rho'} \left(A_-(n,\hbar)\bullet-\rho_0'\right)\right)^*\right\|_{\mathcal{L}(L^(\mathbb{R}^d))},
\end{equation}
where $M$ is a constant depending only on $Q$. We underline that we can choose $\rho$ and $\rho'$ varying in $[-M,M]^{2d}$ (with $M$ fixed) and it will be important to have an uniform bound. 

\subsubsection{Study of the evolution for positive times}
To study the norm of the previous operator, we will first rewrite the operator $\Op_{1}^w\left(\tilde{\Gamma}_{\rho}\left(A_+(n,\hbar)\bullet-\rho_0\right)\right)$ under a more compact form. To do this, define now the Fourier transform of $\tilde{\Gamma}_{\rho}(A_+(n,\hbar)\bullet-\rho_0)$ along the impulsion variable, i.e. for $\rho:=(\rho^1,\rho^2)\in[-M,M]^{2d}$ and $\rho_0:=(\rho_0^1,\rho_0^2)\in\mathbb{R}^{2d}$,
$$
\Gamma_{\rho,\rho_0}^{n,+}(x,\xi):=\frac{1}{(2\pi)^d}\int_{\mathbb{R}^{d}}\tilde{\Gamma}_{\rho}\left(A_+(n,\hbar)\left(\begin{array}{c}x\\ \eta\end{array}\right)-\rho_0\right)e^{\imath\langle \xi|\eta\rangle}d\eta.
$$
With this notation, we can rewrite
$$\forall u\in L^2(\mathbb{R}^d),\ \Op_{1}^w\left(\tilde{\Gamma}_{\rho}\left(A_+(n,\hbar)\bullet-\rho_0\right)\right)u(x)=\int_{\mathbb{R}^d}\Gamma_{\rho,\rho_0}^{n,+}\left(\frac{x+y}{2},x-y\right)u(y)dy.$$
For future purpose, we need to have a precise estimate on the kernel of this operator. Using the Oseledets decomposition of section~\ref{s:oseledets-sympl}, we introduce the notation $(x,\xi):=(\tilde{x}_0,\cdots,\tilde{x}_r,\tilde{\xi}_0,\cdots,\tilde{\xi}_r)\in\mathbb{R}^{2d}$ where $(\tilde{x}_i,\tilde{\xi}_i)$ is an element of $\mathbb{R}^{2d_i}$. Recall also that the matrix $A_i$ that appears in the $\diamond$-decomposition of $A$ is of the form $\text{diag}(D_i,D_i^{-1*})$ for $1\leq i\leq r$. Define now, for $1\leq i\leq r$,
$$\tilde{A}_{i,+}(n,\hbar):=\text{diag}\left(\hbar^{\frac{\lambda_{\max}-\lambda_i^+}{2\lambda_{\max}}}D_i^{\frac{n}{2}},\hbar^{-\frac{\lambda_{\max}-\lambda_i^+}{2\lambda_{\max}}}D_i^{\frac{n}{2}}\right)$$
and
$$\tilde{A}_{0,+}(n,\hbar):=\text{diag}\left(\hbar^{\frac{1}{2}}\text{Id}_{\mathbb{R}^{d_0}},\hbar^{-\frac{1}{2}}\text{Id}_{\mathbb{R}^{d_0}}\right).$$
These matrices allows us to bound the kernel of the operator and to precise where it is large (or not):
\begin{lemm}\label{l:seminorm} Let $L$ be a positive integer. There exists a constant $C_L>0$ such that for every $\rho$ in $[-M,M]^{2d}$ and every $\rho_0:=(\tilde{\rho}_0^{1,0},\cdots,\tilde{\rho}_0^{1,r},\tilde{\rho}_0^{2,0},\cdots,\tilde{\rho}_0^{2,r})\in\mathbb{R}^{2d}$, one has, for every $0\leq n\leq  (1-\epsilon)|\log\hbar|/\lambda_{\max}$ and every $ (x,\xi)\in\mathbb{R}^{2d}$,
$$\prod_{i=0}^r\left(1+\left\|\tilde{A}_{i,+}(n,\hbar)\left(\begin{array}{c}\tilde{x}_i-\tilde{\rho}_0^{1,i}\\ \tilde{\xi}_i\end{array}\right)\right\|^2\right)^L\left|\Gamma_{\rho,\rho_0}^{n,+}(x,\xi)\right|\leq C_L\hbar^{-\frac{d}{2}}\exp\left(\frac{1}{2}\left(n+\frac{\log\hbar}{\lambda_{\max}}\right)\Lambda_+\right),$$
where $\Lambda_+:=\sum_{i=1}^rd_i\lambda_i^+.$
\end{lemm}
\emph{Proof.} For variables corresponding to $1\leq i\leq r$, the estimate follows from a direct computation and uses the fact that $\rho$ varies in a fixed compact set. We have to be more careful in the case $i=0$. For every $\epsilon'>0$, there exists a constant $C_{\epsilon'}>0$ such that, for every $(\tilde{x}_0,\tilde{\xi}_0)$ in $\mathbb{R}^{2d_0}$, one has
$$\forall n\geq0,\ C_{\epsilon'}^{-1}e^{-n\epsilon'}\|(\tilde{x}_0,\tilde{\xi}_0)\|\leq\left\|A_0^{\frac{n}{2}}\left(\begin{array}{c}\tilde{x}_0\\ \tilde{\xi}_0\end{array}\right)\right\|\leq C_{\epsilon'}^{}e^{n\epsilon'}\|(\tilde{x}_0,\tilde{\xi}_0)\|.$$
These estimates allow us to obtain the bounds we need (using also the fact that we restrict ourselves to times $n$ that are at most logarithmic in $\hbar$).$\square$

\subsubsection{Study of the evolution for negative times} In the previous paragraph, we studied the norm of the operator for positive times. We can also rewrite 
$$\forall u\in L^2(\mathbb{R}^d),\ \Op_{1}^w\left(\tilde{\Gamma}_{\rho'}\left(A_-(n,\hbar)\bullet-\rho_0'\right)\right)u(x)=\int_{\mathbb{R}^d}\overline{\Gamma}_{\rho',\rho_0'}^{-n,+}\left(\frac{x+y}{2},y-x\right)u(y)dy.$$
If we define $\tilde{A}_{i,-}(n,\hbar):=\tilde{A}_{i,+}(-n,\hbar)$ for $0\leq i\leq r$, we also have the following lemma:
\begin{lemm}\label{l:seminormneg} Let $L$ be a positive integer. There exists a constant $C_L>0$ such that for every $\rho'$ in $[-M,M]^{2d}$ and every $\rho_0:=(\tilde{\rho}_0^{1,0},\cdots,\tilde{\rho}_0^{1,r},\tilde{\rho}_0^{2,0},\cdots,\tilde{\rho}_0^{2,r})\in\mathbb{R}^{2d}$, one has, for every $0\leq n\leq (1-\epsilon)|\log\hbar|/\lambda_{\max}$ and every $ (x,\xi)\in\mathbb{R}^{2d}$,
$$\prod_{i=0}^r\left(1+\left\|\tilde{A}_{i,-}(n,\hbar)\left(\begin{array}{c}\tilde{x}_i-\tilde{\rho}_0^{1,i}\\ \tilde{\xi}_i\end{array}\right)\right\|^2\right)^L\left|\Gamma_{\rho,\rho_0}^{-n,+}(x,\xi)\right|\leq C_L\hbar^{-\frac{d}{2}}\exp\left(\frac{1}{2}\left(\frac{\log\hbar}{\lambda_{\max}}-n\right)\Lambda_+\right),$$
where $\Lambda_+:=\sum_{i=1}^rd_i\lambda_i^+.$
\end{lemm}

\subsubsection{Estimate of the norm}
In order to compute a bound on the norm of the operator, we consider two elements $\phi_1$ and $\phi_2$ in $L^2(\mathbb{R}^d)$ and we want to estimate
$$C_{\phi_1,\phi_2}(n):=\left\langle\phi_1,\Op_{1}^w\left(\tilde{\Gamma}_{\rho}\left( A_+(n,\hbar)\bullet-\rho_0\right)\right)\Op_{1}^w\left(\tilde{\Gamma}_{\rho'}\left( A_-(n,\hbar)\bullet-\rho_0'\right)\right)^*\phi_2\right\rangle_{L^{2}(\mathbb{R}^d)}.$$
With the notations of the previous paragraphs, one has
$$C_{\phi_1,\phi_2}(n):=\int_{\mathbb{R}^{3d}}\Gamma_{\rho,\rho_0}^{n,+}\left(\frac{x+y}{2},x-y\right)\overline{\Gamma}_{\rho',\rho_0'}^{-n,+}\left(\frac{z+y}{2},y-z\right)\overline{\phi}_1(x)\phi_2(z)dxdydz.$$
Fix now some positive (small) number $\delta$ and introduce the two following subsets of $\mathbb{R}^3$:
$$X_{\delta}^+(n):=\left\{(x,y,z):\forall 0\leq i\leq r \left\|\tilde{A}_{i,+}(n,\hbar)\left(\begin{array}{c}\frac{\tilde{x}_i+\tilde{y}_i}{2}-\tilde{\rho}_0^{1,i}\\ \tilde{x}_i-\tilde{y}_i\end{array}\right)\right\|^2\leq\hbar^{-\delta}\right\}$$
and
$$X_{\delta}^-(n):=\left\{(x,y,z):\forall 0\leq i\leq r \left\|\tilde{A}_{i,-}(n,\hbar)\left(\begin{array}{c}\frac{\tilde{z}_i+\tilde{y}_i}{2}-\tilde{\rho}_0^{'1,i}\\ \tilde{y}_i-\tilde{z}_i\end{array}\right)\right\|^2\leq\hbar^{-\delta}\right\}.$$
We also define $X_{\delta}(n):=X_{\delta}^+(n)\cap X_{\delta}^-(n)$. Thanks to lemmas~\ref{l:seminorm} and~\ref{l:seminormneg}, we know that outside the set $X_{\delta}(n)$, the kernel of the operator is small in $\hbar$. Precisely, we can prove the following lemma:
\begin{lemm} Let $\delta$ be a (small) positive real number. Let $\phi_1$ and $\phi_2$ be two elements in $L^2(\mathbb{R}^d)$. One has, for $0\leq n\leq (1-\epsilon)|\log\hbar|/\lambda_{\max}$,
$$\int_{\mathbb{R}^{3d}\backslash X_{\delta}(n)}\Gamma_{\rho,\rho_0}^{n,+}\left(\frac{x+y}{2},x-y\right)\overline{\Gamma}_{\rho',\rho_0'}^{-n,+}\left(\frac{z+y}{2},y-z\right)\overline{\phi}_1(x)\phi_2(z)dxdydz=\mathcal{O}(\hbar^{\infty})\|\phi_1\|_{L^2(\mathbb{R}^d)}\|\phi_2\|_{L^2(\mathbb{R}^d)},$$
where the remainder is uniform for $\rho,\rho'\in[-M,M]^{2d}$ and $\rho_0,\rho_0'\in\mathbb{R}^{2d}$.
\end{lemm}
\emph{Proof.} Thanks to the Cauchy Schwarz inequality, it is sufficient to give an estimate on 
$$\int_{\mathbb{R}^{3d}\backslash X_{\delta}(n)}\left|\Gamma_{\rho,\rho_0}^{n,+}\left(\frac{x+y}{2},x-y\right)\overline{\Gamma}_{\rho',\rho_0'}^{-n,+}\left(\frac{z+y}{2},y-z\right)\right|\left|\overline{\phi}_1(x)\right|^2dxdydz.$$
According to lemmas~\ref{l:seminorm} and~\ref{l:seminormneg}, we know that for every integer $L$, there exists a constant $C_L>0$ such that, for every $0\leq n\leq|\log\hbar|/\lambda_{\max}$ and every $(x,y,z)$ in $\mathbb{R}^{3d}$
$$\left|\Gamma_{\rho,\rho_0}^{n,+}\left(\frac{x+y}{2},x-y\right)\overline{\Gamma}_{\rho',\rho_0'}^{-n,+}\left(\frac{z+y}{2},y-z\right)\right|$$
$$\leq C_L\hbar^{-d}\prod_{i=0}^r\left(1+\left\|\tilde{A}_{i,+}(n,\hbar)\left(\begin{array}{c}\frac{\tilde{x}_i+\tilde{y}_i}{2}-\tilde{\rho}_0^{1,i}\\ \tilde{x}_i-\tilde{y}_i\end{array}\right)\right\|^2\right)^{-L}\prod_{i=0}^r\left(1+\left\|\tilde{A}_{i,-}(n,\hbar)\left(\begin{array}{c}\frac{\tilde{z}_i+\tilde{y}_i}{2}-\tilde{\rho}_0^{'1,i}\\ \tilde{y}_i-\tilde{z}_i\end{array}\right)\right\|^2\right)^{-L}.$$
Under the extra assumption that $(x,y,z)$ in $\mathbb{R}^{3d}\backslash X_{\delta}$, one knows that
$$\left|\Gamma_{\rho,\rho_0}^{n,+}\left(\frac{x+y}{2},x-y\right)\overline{\Gamma}_{\rho',\rho_0'}^{-n,+}\left(\frac{z+y}{2},y-z\right)\right|$$
$$\leq C_L\hbar^{\delta L-d}\prod_{i=0}^r\left(1+\left\|\tilde{A}_{i,+}(n,\hbar)\left(\begin{array}{c}\frac{\tilde{x}_i+\tilde{y}_i}{2}-\tilde{\rho}_0^{1,i}\\ \tilde{x}_i-\tilde{y}_i\end{array}\right)\right\|^2\right)^{-d}\prod_{i=0}^r\left(1+\left\|\tilde{A}_{i,-}(n,\hbar)\left(\begin{array}{c}\frac{\tilde{z}_i+\tilde{y}_i}{2}-\tilde{\rho}_0^{'1,i}\\ \tilde{y}_i-\tilde{z}_i\end{array}\right)\right\|^2\right)^{-d}.$$
In the allowed range of times $n$ one can check that, for some uniform constant $D$, one has 
$$\int_{\mathbb{R}^{2d}}\prod_{i=0}^r\left(1+\left\|\tilde{A}_{i,-}(n,\hbar)\left(\begin{array}{c}\frac{\tilde{z}_i+\tilde{y}_i}{2}-\tilde{\rho}_0^{'1,i}\\ \tilde{y}_i-\tilde{z}_i\end{array}\right)\right\|^2\right)^{-d}dydz=\mathcal{O}(\hbar^{-D}).$$
Combining these last two estimates, we find that, for every $L>0$,
$$\int_{\mathbb{R}^{3d}\backslash X_{\delta}(n)}\left|\Gamma_{\rho,\rho_0}^{n,+}\left(\frac{x+y}{2},x-y\right)\overline{\Gamma}_{\rho',\rho_0'}^{-n,+}\left(\frac{z+y}{2},y-z\right)\right|\left|\overline{\phi}_1(x)\right|^2dxdydz=\mathcal{O}(\hbar^{\delta L-D-d}).\square$$

\subsubsection{The conclusion} According to the previous paragraph, we know that, modulo a remainder of order $\mathcal{O}(\hbar^{\infty})\|\phi_1\|_{L^2(\mathbb{R}^d)}\|\phi_2\|_{L^2(\mathbb{R}^d)}$, the quantity $C_{\phi_1,\phi_2}(n)$ is equal to
\begin{equation}\label{e:maincontribution}\int_{X_{\delta}(n)}\Gamma_{\rho,\rho_0}^{n,+}\left(\frac{x+y}{2},x-y\right)\overline{\Gamma}_{\rho',\rho_0'}^{-n,+}\left(\frac{z+y}{2},y-z\right)\overline{\phi}_1(x)\phi_2(z)dxdydz.\end{equation}
For the sake of simplicity, we now restrict ourselves to the case\footnote{The reader can check that the bound is optimal in this case.} $n=n_E(\hbar)=[(1-\epsilon)|\log\hbar|]$. According to lemmas~\ref{l:seminorm} and~\ref{l:seminormneg}, one knows that there exists a constant $C>0$ such that
$$\left|(\ref{e:maincontribution})\right|\leq C\hbar^{-d+\frac{\Lambda_+}{\lambda_{\max}}}\int_{X_{\delta}(n_{E}(\hbar))}|\phi_1(x)\phi_2(z)|dxdydz.$$
We have then
\begin{equation}\label{e:laststep-lead-term}\left|(\ref{e:maincontribution})\right|\leq C\hbar^{-d+\frac{\Lambda_+}{\lambda_{\max}}}\left(\int_{X_{\delta}(n_{E}(\hbar))}|\phi_1(x)|^2dxdydz\right)^{\frac{1}{2}}\left(\int_{X_{\delta}(n_{E}(\hbar))}|\phi_2(z)|^2dxdydz\right)^{\frac{1}{2}},\end{equation}
We will estimate these two integrals and to do this, we will distinguish two cases:
\begin{itemize}
\item the variables such that $2\lambda_i^+-\lambda_{\max}>0$;
\item the variables such that $2\lambda_i^+-\lambda_{\max}\leq0$.
\end{itemize}
\textbf{The case $2\lambda_i^+-\lambda_{\max}>0$.} Define, for $i$ such that $2\lambda_i^+-\lambda_{\max}>0$ and $\tilde{x}_i$ in $\mathbb{R}^{d_i}$, the sets
$$X_{\delta}^{+}(\tilde{x}_i,n_E(\hbar)):=\left\{(\tilde{y}_i,\tilde{z}_i): \left\|\tilde{A}_{i,+}(n_E(\hbar),\hbar)\left(\begin{array}{c}\frac{\tilde{x}_i+\tilde{y}_i}{2}-\tilde{\rho}_0^{1,i}\\ \tilde{x}_i-\tilde{y}_i\end{array}\right)\right\|^2\leq\hbar^{-\delta}\right\}$$
and
$$X_{\delta}^{-}(\tilde{x}_i,n_E(\hbar)):=\left\{(\tilde{y}_i,\tilde{z}_i): \left\|\tilde{A}_{i,-}(n_E(\hbar),\hbar)\left(\begin{array}{c}\frac{\tilde{z}_i+\tilde{y}_i}{2}-\tilde{\rho}_0^{'1,i}\\ \tilde{y}_i-\tilde{z}_i\end{array}\right)\right\|^2\leq\hbar^{-\delta}\right\}.$$
For every $\tilde{x}_i$ in $\mathbb{R}^{d_i}$, we can verify that the optimal bound on the volume is given by
$$\text{Vol}\left(X_{\delta}^{+}(\tilde{x}_i,n_E(\hbar))\cap X_{\delta}^{-}(\tilde{x}_i,n_E(\hbar))\right)\leq \tilde{C} \hbar^{-d_i\delta}\hbar^{d_i\left(1-\frac{\lambda_i^+}{\lambda_{\max}}\right)},$$
where $\tilde{C}$ is some uniform constant. If we do the same thing but exchange the roles played by $\tilde{x}_i$ and $\tilde{z}_i$, we can introduce the sets
$$X_{\delta}^{+}(\tilde{z}_i,n_E(\hbar)):=\left\{(\tilde{x}_i,\tilde{y}_i): \left\|\tilde{A}_{i,+}(n_E(\hbar),\hbar)\left(\begin{array}{c}\frac{\tilde{x}_i+\tilde{y}_i}{2}-\tilde{\rho}_0^{1,i}\\ \tilde{x}_i-\tilde{y}_i\end{array}\right)\right\|^2\leq\hbar^{-\delta}\right\}$$
and
$$X_{\delta}^{-}(\tilde{z}_i,n_E(\hbar)):=\left\{(\tilde{x}_i,\tilde{y}_i): \left\|\tilde{A}_{i,-}(n_E(\hbar),\hbar)\left(\begin{array}{c}\frac{\tilde{z}_i+\tilde{y}_i}{2}-\tilde{\rho}_0^{'1,i}\\ \tilde{y}_i-\tilde{z}_i\end{array}\right)\right\|^2\leq\hbar^{-\delta}\right\}.$$
We verify that the optimal bound on the volume is given by
$$\text{Vol}\left(X_{\delta}^{+}(\tilde{z}_i,n_E(\hbar))\cap X_{\delta}^{-}(\tilde{z}_i,n_E(\hbar))\right)\leq \tilde{C} \hbar^{-\frac{d_i\delta}{2}}|\det \tilde{A}_{i,+}(n_E(\hbar),\hbar)|^{-1}=\tilde{C} \hbar^{-d_i\delta}e^{-n_E(\hbar)d_i\lambda_i^+}.$$
For this bound we used the fact that the Lyapunov exponent satisfy $2\lambda_i^+-\lambda_{\max}$: the volume $\text{Vol}\left(X_{\delta}^{+}(\tilde{z}_i,n_E(\hbar))\right)$ is already bounded by this quantity and it is the optimal bound we can get (regarding the fact that $2\lambda_i^+-\lambda_{\max}$). These estimates will allow us to treat the variables corresponding to indices $i$ such that $2\lambda_i^+-\lambda_{\max}>0$ in the right hand side of~(\ref{e:laststep-lead-term}).\\
\\
\textbf{The case $2\lambda_i^+-\lambda_{\max}\leq 0$.} We now treat the case of the other variables. We fix such a $i$. We also use the same auxiliary sets for $\tilde{x}_i$ and $\tilde{z}_i$ in $\mathbb{R}^{d_i}$. For every $\tilde{x}_i$ in $\mathbb{R}^{d_i}$, we can verify that, as in the previous case, the optimal bound on the volume is given by
$$\text{Vol}\left(X_{\delta}^{+}(\tilde{x}_i,n_E(\hbar))\cap X_{\delta}^{-}(\tilde{x}_i,n_E(\hbar))\right)\leq \tilde{C} \hbar^{-d_i\delta}\hbar^{d_i\left(1-\frac{\lambda_i^+}{\lambda_{\max}}\right)},$$
where $\tilde{C}$ is some uniform constant. The difference with the previous case is that, as $2\lambda_i^+-\lambda_{\max}\leq 0$, we can not obtain a better bound in the case of $\tilde{z}_i$. It means that we have, for every $\tilde{z}_i$, the optimal bound on the volume is given by
$$\text{Vol}\left(X_{\delta}^{+}(\tilde{z}_i,n_E(\hbar))\cap X_{\delta}^{-}(\tilde{z}_i,n_E(\hbar))\right)\leq \tilde{C} \hbar^{-d_i\delta}\hbar^{d_i\left(1-\frac{\lambda_i^+}{\lambda_{\max}}\right)}.$$
\textbf{Combining the different estimates.} Using the previous definitions, we have the following inequality
$$\int_{\tilde{X}_{\delta}(n_{E}(\hbar))}|\phi_1(x)|^2dxdydz\leq\int_{\tilde{X}_{\delta}(n_{E}(\hbar))}|\phi_1(x)|^2\prod_{i=0}^r\text{Vol}\left(X_{\delta}^{+}(\tilde{x}_i,n_E(\hbar))\cap X_{\delta}^{-}(\tilde{x}_i,n_E(\hbar))\right)d\tilde{x}_0\cdots d\tilde{x}_r.$$
With our previous estimates, we find that
$$\int_{\tilde{X}_{\delta}(n_{E}(\hbar))}|\phi_1(x)|^2dxdydz\leq\tilde{C}^2\hbar^{-d\delta+d-\frac{\Lambda_+}{\lambda_{\max}}}\|\phi_1\|_{L^2(\mathbb{R}^d)}^2.$$
where we recall that $\Lambda_+:=\sum_{i=1}^r d_i\lambda_i^+$. For the other integral, we find that
$$\int_{\tilde{X}_{\delta}(n_{E}(\hbar))}|\phi_2(z)|^2dxdydz\leq\tilde{C}^2\hbar^{-d\delta}\hbar^{\sum_{i:2\lambda_i^+-\lambda_{\max}\leq 0}d_i\left(1-\frac{\lambda_i^+}{\lambda_{\max}}\right)}e^{-n_E(\hbar)\sum_{i:2\lambda_i^+-\lambda_{\max}>0}d_i\lambda_i^+}\|\phi_2\|_{L^2(\mathbb{R}^d)}^2.$$
Finally, using~(\ref{e:laststep-lead-term}), we find that, for $n=n_E(\hbar)$,
$$\left|(\ref{e:maincontribution})\right|\leq C\tilde{C}^2\hbar^{-d\delta-\epsilon\frac{\Lambda_+}{\lambda_{\max}}}\exp\left(\sum_{i:2\lambda_i^+-\lambda_{\max}>0}d_i\left(\lambda_i^+-\frac{\lambda_{\max}}{2}\right)\frac{\log\hbar}{\lambda_{\max}}\right)\|\phi_1\|_{L^2(\mathbb{R}^d)}\|\phi_2\|_{L^2(\mathbb{R}^d)}.\square$$

\subsection{Negligible terms}\label{s:negligible} In the previous section, we have estimated the term that is supposed to be the leading term in the operator norm. We will prove in this section that some of the other terms are negligible and we will conclude in the next section using the Coltar Stein theorem~\cite{DS}.\\
As $A$ is an element in $SL(2d,\mathbb{Z})$, we can consider $T_{\rho}\left(G_{\hbar}\circ A^{\frac{n}{2}}\right)^{\rho_0}$ instead of $T_{\rho}\left(G_{\hbar}^{\rho_0}\right)\circ A^{\frac{n}{2}}$ and $T_{\rho'}\left(G_{\hbar}\circ A^{-\frac{n}{2}}\right)^{\rho_0'}$ instead of  $T_{\rho'}\left(G_{\hbar}^{\rho_0'}\right)\circ A^{-\frac{n}{2}}$ where $\rho$, $\rho'$, $\rho_0$ and $\rho_0'$ vary in $\mathbb{T}^{2d}$ (see the expression~(\ref{e:c(A,n)}) we want to estimate). The observables $T_{\rho}\left(G_{\hbar}\circ A^{\frac{n}{2}}\right)^{\rho_0}$ and $T_{\rho'}\left(G_{\hbar}\circ A^{-\frac{n}{2}}\right)^{\rho_0'}$ are defined as infinite sums over $r$ in $\mathbb{Z}^{2d}$. We will now show that for $r\neq r'$, the norm of the product operator is a $\mathcal{O}(\hbar^{\infty})$. First, for simplicity of notations, we introduce the following notations,
 for $w\in\mathbb{R}^{2d}$,
$$F_+^{(n)}(w):=G\left(B(\hbar)A^{\frac{n}{2}}\left(w-\rho_0-\pi\hbar\rho\right)\right)e^{2\imath\pi\langle A^{\frac{n}{2}}w|\rho\rangle}$$
and
$$F_-^{(n)}(w):=G\left(B(\hbar)A^{-\frac{n}{2}}\left(w-\pi\hbar\rho'\right)\right)e^{-2\imath\pi\langle A^{-\frac{n}{2}}w|\rho'\rangle}.$$
We underline that we have taken $\rho_0'=0$ without loss of generality (see the expression~(\ref{e:c(A,n)}) we want to estimate). Moreover, we can also suppose that $\rho_0$ is an element in $[-1/2,1/2]^{2d}$. We now estimate the norm of two translated operators with $r\neq r'$. To do this, we write the exact formula for the Moyal product (see~\cite{EZ}-chapter $4$), for $r$ and $r'$ in $\mathbb{Z}^{2d}$,
$$A_{r,r'}(w):=F_+^{(n)}(\bullet+r)\sharp F_-^{(n)}(\bullet+r')(w)=\int_{\mathbb{R}^{4d}}F_+^{(n)}(w+w_1+r)F_-^{(n)}(w+w_2+r')e^{-\frac{2\imath}{\hbar}\langle w_1,Jw_2\rangle}\frac{dw_1dw_2}{(\pi\hbar)^{2d}}.$$
Let $\chi(w_1,w_2)$ be a smooth function on $\mathbb{R}^{4d}$ compactly supported in a small neighborhood of $0$. We fix some small positive number $\epsilon'$ and we suppose that $\chi$ is equal to $1$ on the set $\{\|w_1\|_2\leq\epsilon'\ \text{and}\ \|w_2\|_2\leq\epsilon'\}$ and to $0$ outside $\{\|w_1\|_2\leq2\epsilon'\ \text{and}\ \|w_2\|_2\leq2\epsilon'\}$. Using this cutoff, we can split the integral in two parts
$$A_{r,r'}^1(w):=\int_{\mathbb{R}^{4d}}\chi(w_1,w_2)F_+^{(n)}(w+w_1+r)F_-^{(n)}(w+w_2+r')e^{-\frac{2\imath}{\hbar}\langle w_1,Jw_2\rangle}\frac{dw_1dw_2}{(\pi\hbar)^{2d}}$$
and
$$A_{r,r'}^2(w):=\int_{\mathbb{R}^{4d}}(1-\chi(w_1,w_2))F_+^{(n)}(w+w_1+r)F_-^{(n)}(w+w_2+r')e^{-\frac{2\imath}{\hbar}\langle w_1,Jw_2\rangle}\frac{dw_1dw_2}{(\pi\hbar)^{2d}}.$$
We will now prove that these two symbols are in the class $S^{-\infty}(1)$ with an explicit control on the norm of the derivatives depending on $r$ and $r'$.

\subsubsection{Class of $A_{r,r'}^2$} We know that the integral defining $A_{r,r'}^2$ is over variables $(w_1,w_2)$ that satisfy
$$\left\|w_1\right\|_2>\epsilon'\ \text{or}\ \left\|w_2\right\|_2>\epsilon'.$$
Thanks to this last property, we are able to use the (non)-stationary phase property. To do this, we introduce the operators $$L:=\frac{\hbar}{2\imath}\left\langle\frac{ w_1}{\left\|w_1\right\|_2^2},Jd_{w_2}\right\rangle\ \text{or}\ L':=-\frac{\hbar}{2\imath}\left\langle\frac{ Jw_2}{\left\|w_2\right\|_2^2},d_{w_1}\right\rangle.$$
Using the fact that $L(e^{-\frac{2\imath}{\hbar}\langle w_1,Jw_2\rangle})=L'(e^{-\frac{2\imath}{\hbar}\langle w_1,Jw_2\rangle})=e^{-\frac{2\imath}{\hbar}\langle w_1,Jw_2\rangle}$ and performing integration by parts, we find that the observable $A_{r,r'}^2(w)$ is a $\mathcal{O}(\hbar^{\infty})$ as long as $0\leq n\leq \frac{1-\epsilon}{\lambda_{\max}}|\log\hbar|$ (the derivatives of $F_+^{(n)}$ and $F_-^{(n)}$ are bounded by some $\mathcal{O}(\hbar^{-1+\frac{\epsilon}{2}})$ for this range of times). Moreover, we can make other integrations by parts using the operators
$$L_r:=\frac{1+\frac{\hbar}{2\imath}\left\langle w+r,Jd_{w_2}\right\rangle}{1+\|w+r\|_2^2}\ \text{and}\ L_{r'}':=\frac{1-\frac{\hbar}{2\imath}\left\langle w+r',d_{w_1}\right\rangle}{1+\|w+r'\|_2^2}.$$
We verify then that, for every $M$ in $\mathbb{N}$, there exists a constant $C_M$ such that
$$\forall r\neq r'\in\mathbb{Z}^{2d},\ \forall w\in \mathbb{R}^{2d},\ |A_{r,r'}^2(w)|\leq \frac{C_M\hbar^M}{(1+\|w+r'\|_2)^{2d}(1+\|w+r\|_2)^{2d}}.$$
Making the same computations, we find the same properties hold for any derivative of $A_{r,r'}^2$. In particular, we know that the symbol $\sum_{r\neq r'}A_{r,r'}^2$ is in the class $S^{-\infty}(1)$, as long as $0\leq n\leq \frac{1-\epsilon}{\lambda_{\max}}|\log\hbar|$.

\subsubsection{Class of $A_{r,r'}^1$} For $\hbar$ small enough and for $w_1$ on the support of $\chi$, we know that the observable $F_+^{(n)}(w+w_1+r)$ is gaussian and centered on a point in the ball $B(r,3\epsilon'+1/2)$. Moreover, for $w_2$ on the support of $\chi$, the other observable $F_-^{(n)}(w+w_2+r')$ is gaussian and centered on a point in the ball $B(r',3\epsilon')$ (again when $\hbar$ is small enough). As we made the assumption that $r\neq r'$, we also know that $\|r-r'\|_2\geq 1$. If we restrict ouselves to $0\leq n\leq (1-\epsilon)|\log\hbar|/\lambda_{\max}$, the variance of the two gaussian observables is of order at most $\mathcal{O}(\hbar^{\epsilon}).$ These different observations tell us that the observable $F_+^{(n)}(w+w_1+r)$ is exponentially small in $\hbar$ when $F_-^{(n)}(w+w_2+r')$ is large. The converse is also true. In particular, we know that $|A_{r,r'}^1(w)|=\mathcal{O}(\hbar^{\infty})$ (uniformly for $w$ in $\mathbb{R}^{2d}$). In fact, we can even be more precise and we can verify that, for every $L>0$, 
$$(1+\|w+r'\|_2)^{2d}(1+\|w+r\|_2)^{2d}|A_{r,r'}^1(w)|=\mathcal{O}(\hbar^L),$$
where the constant involved is uniform for $r\neq r'$ in $\mathbb{Z}^{2d}$, $w$ in $\mathbb{R}^{2d}$ and $0\leq n\leq(1-\epsilon)|\log\hbar|/\lambda_{\max}$. Finally, we underline that the same method allows to derive the same on the derivatives of $A_{r,r'}^1$. In particular, the symbol $\sum_{r\neq r'} A_{r,r'}^1$ is in the class $S^{-\infty}(1)$.

\subsubsection{Applying Calder\'on Vaillancourt theorem} Using the two previous paragraphs, we know that the symbol $\sum_{r\neq r'}A_{r,r'}$ is in the class $S^{-\infty}(1)$. Thanks to the Calder\'on Vaillancourt theorem (see equation~(\ref{e:calderon-vail})), we know that, as long as $n\leq \frac{1-\epsilon}{\lambda_{\max}}|\log\hbar|$,
$$\left\|\Op_{\hbar}\left(\sum_{r\neq r'}A_{r,r'}\right)\right\|_{\mathcal{L}(L^2(\mathbb{R}^d))}=\mathcal{O}(\hbar^{\infty}).$$
Finally, we can derive that
$$\left\|\Op_{\hbar}^w(T_{\rho}(G_{\hbar}\circ A^{\frac{n}{2}})^{\rho_0})\Op_{\hbar}^w(T_{\rho'}(G_{\hbar}\circ A^{-\frac{n}{2}})^{\rho_0'})^*\right\|_{\mathcal{L}(L^2(\mathbb{R}^{d}))}\ \ \ \ \ \ \ \ \ \ \ \ \ \ \ \ \ \ \ \ \ \ \ \ \ \ \ \ \ \ \ \ $$
$$\ \ \ \ \ \ \ \ \ \ \ \ \ \ \ \ \ \ \ \ \ \ \ \ \ \ \ \ \ \ =2^d|\det B(\hbar)|\left\|\Op_{\hbar}^w\left(T_{0}\left(F_+^{(n)}\sharp F_-^{(n)}\right)\right)\right\|_{\mathcal{L}(L^2(\mathbb{R}^{d}))}+O(\hbar^{\infty}).$$

\subsection{Applying Coltar-Stein theorem}\label{s:Coltar}
 To summarize, we have shown that in order to prove theorem~\ref{t:main-est}, we only need to get an estimate on the norm of the operator
$$\Op_{\hbar}^w(T_{0}(F_+^{(n)}\sharp F_-^{(n)}))=\sum_{r\in\mathbb{Z}^{2d}}U_{\hbar}(r)\Op_{\hbar}^w(F_+^{(n)}\sharp F_-^{(n)})U_{\hbar}(r)^*,$$
where we used the notations of the previous section. Moreover, proposition~\ref{p:leadterm} shows that the norm of $\Op_{\hbar}^w(F_+\sharp F_-)$ is bounded by the expected quantity. It remains to show that these two properties are sufficient to prove the main theorem. To do this, we define $\mathbf{A}_r:=U_{\hbar}(r)\Op_{\hbar}^w(F_+^{(n)}\sharp F_-^{(n)})U_{\hbar}(r)^*$. Our goal is to give a bound on the two following quantities:
$$\sup_r\sum_{r'\in\mathbb{Z}^{2d}}\|\mathbf{A}_r^*\mathbf{A}_{r'}\|_{\mathcal{L}(L^2(\mathbb{R}^d))}^{\frac{1}{2}}\ \text{and}\ \sup_r\sum_{r'\in\mathbb{Z}^{2d}}\|\mathbf{A}_r\mathbf{A}_{r'}^*\|_{\mathcal{L}(L^2(\mathbb{R}^d))}^{\frac{1}{2}}.$$
If we are able to prove that both quantities are bounded by the same quantity $C$, Coltar-Stein theorem will tell us that $C$ is a bound on the norm of $\mathbf{A}:=\sum_{r\in\mathbb{Z}^{2d}}\mathbf{A}_r$. Regarding this goal, we write
$$\mathbf{A}_r^*\mathbf{A}_{r'}=\Op_{\hbar}^w((F_-^{(n)})^{r})^*\Op_{\hbar}^w((F_+^{(n)})^{r})^*\Op_{\hbar}^w((F_+^{(n)})^{r'})\Op_{\hbar}^w((F_-^{(n)})^{r'}).$$
Using the exact Egorov property, we know that $\displaystyle\left\|\Op_{\hbar}^w((F_-^{(n)})^r)\right\|_{\mathcal{L}(L^2(\mathbb{R}^d))}=\left\|\Op_{\hbar}^w(F_-\circ A^{\frac{n}{2}})\right\|_{\mathcal{L}(L^2(\mathbb{R}^d))}$. As long as $|n|\leq \frac{1-\epsilon}{\lambda_{\max}}|\log \hbar|$, the symbol $F_-\circ A^{\frac{n}{2}}$ remains in the class of symbol $S^0_{\frac{1}{2}}(1)$ and the involved constants are uniform for $n$ in the allowed interval. In particular, we know that there exists an uniform constant $C$ (independent of $r$) such that $\displaystyle\left\|\Op_{\hbar}^w((F_-^{(n)})^r)\right\|_{\mathcal{L}(L^2(\mathbb{R}^d))}\leq C$. In particular, we have that $\displaystyle\|\mathbf{A}_r^*\mathbf{A}_{r'}\|_{\mathcal{L}(L^2(\mathbb{R}^d))}^{\frac{1}{2}}\leq C \left\|\Op_{\hbar}^w((F_+^{(n)})^r)^*\Op_{\hbar}^w((F_+^{(n)})^{r'})\right\|_{\mathcal{L}(L^2(\mathbb{R}^d))}^{\frac{1}{2}}.$ We will use this estimate to bound the sum over $r\neq r'$. Using the same method as in section~\ref{s:negligible}, we find that for every $M$ in $\mathbb{N}$ and every multiindex $\alpha$ in $\mathbb{N}^{2d}$, there exists a constant $C_{M,\alpha}$ such that
$$\forall\ r\neq r',\ \forall w\in\mathbb{R}^{2d},\ |\partial^{\alpha}((F_+^{(n)})^r\sharp (F_+^{(n)})^{r'})(w)|\leq\frac{C_{M,\alpha}\hbar^M}{(1+\|w+r'\|_2)^{2d}(1+\|w+r\|_2)^{2d}}.$$
Finally, according to Calder\'on-Vaillancourt theorem, we know then that, for every $M$ in $\mathbb{N}$, there exists a constant $C_M$, such that
$$\forall\ r\neq r',\ \|\mathbf{A}_r^*\mathbf{A}_{r'}\|_{\mathcal{L}(L^2(\mathbb{R}^d))}^{\frac{1}{2}}\leq C_M\hbar^M (1+\|r-r'\|_2)^{-2d}.$$
In particular, it implies that, for every $r\in\mathbb{Z}^{2d}$,
$$\sum_{r'\in\mathbb{Z}^{2d}}\|\mathbf{A}_r^*\mathbf{A}_{r'}\|_{\mathcal{L}(L^2(\mathbb{R}^d))}^{\frac{1}{2}}=\left\|\Op_{\hbar}^w(F_-^{(n)})^*\Op_{\hbar}^w(F_+^{(n)})^*\Op_{\hbar}^w(F_+^{(n)})\Op_{\hbar}^w(F_-^{(n)})\right\|_{\mathcal{L}(L^2(\mathbb{R}^d))}^{\frac{1}{2}}+O(\hbar^{\infty}).$$
As we have seen with proposition~\ref{p:leadterm} how to to bound the norm $\Op_{\hbar}^w(F_+^{(n)})\Op_{\hbar}^w(F_-^{(n)})$ and as we have $\displaystyle\left\|\Op_{\hbar}^w(F_-^{(n)})^*\Op_{\hbar}^w(F_+^{(n)})^*\right\|_{\mathcal{L}(L^2(\mathbb{R}^d))}=\left\|\Op_{\hbar}^w(F_+^{(n)})\Op_{\hbar}^w(F_-^{(n)})\right\|_{\mathcal{L}(L^2(\mathbb{R}^d))}$, we know that, for $n=n_E(\hbar):=[(1-\epsilon)|\log\hbar|/\lambda_{\max}]$,
$$\sum_{r'\in\mathbb{Z}^{2d}}\|\mathbf{A}_r^*\mathbf{A}_{r'}\|_{\mathcal{L}(L^2(\mathbb{R}^d))}^{\frac{1}{2}}\leq C\hbar^{-\delta-\epsilon\frac{\Lambda_+}{\lambda_{\max}}}\exp\left(-\sum_{i:2\lambda_i^+-\lambda_{\max}>0}d_i\left(\lambda_i^+-\frac{\lambda_{\max}}{2}\right)n_E(\hbar)\right).$$
The same method allows to get the bound $$\sup_r\sum_{r'\in\mathbb{Z}^{2d}}\|\mathbf{A}_r\mathbf{A}_{r'}^*\|_{\mathcal{L}(L^2(\mathbb{R}^d))}^{\frac{1}{2}}\leq C\hbar^{-\delta-\epsilon\frac{\Lambda_+}{\lambda_{\max}}}\exp\left(-\sum_{i:2\lambda_i^+-\lambda_{\max}>0}d_i\left(\lambda_i^+-\frac{\lambda_{\max}}{2}\right)n_E(\hbar)\right).$$ By Coltar-Stein theorem (lemma $7.10$ in~\cite{DS}), we can deduce that, for $n=n_E(\hbar)$,
$$\left\|\Op_{\hbar}^w(T_{0}(F_+^{(n)}\sharp F_-^{(n)}))\right\|_{\mathcal{L}(L^2(\mathbb{R}^d))}\leq C\hbar^{-\delta-\epsilon\frac{\Lambda_+}{\lambda_{\max}}}\exp\left(-\sum_{i:2\lambda_i^+-\lambda_{\max}>0}d_i\left(\lambda_i^+-\frac{\lambda_{\max}}{2}\right)n_E(\hbar)\right).\square$$

\appendix

\section{Proof of lemma~\ref{l:equiv-procedure}}\label{a:eq-proc} In this appendix, we give a proof of lemma~\ref{l:equiv-procedure}. Precisely, we have to verify that the symbols $a$ and $a\star\left(G_{\hbar}\sharp G_{\hbar}\right)$ have the same principal symbol. Using the definition of the Moyal product~\cite{EZ}, we can compute an exact expression of the symbol
$$\begin{array}{ccc}
a\star\left(G_{\hbar}\sharp G_{\hbar}\right)(\rho)&=&\int_{\mathbb{R}^{2d}}a(\rho_0)\int_{\mathbb{R}^{4d}}e^{-\frac{2\imath}{\hbar}\omega(w_1,w_2)}G_{\hbar}(\rho-\rho_0+w_1)G_{\hbar}(\rho-\rho_0+w_2)\frac{dw_1dw_2}{(\pi\hbar)^{2d}}d\rho_0\\
&=&\int_{\mathbb{R}^{2d}}a(\rho_0)\int_{\mathbb{R}^{4d}}e^{-2\imath\pi\omega(w_1,w_2)}G_{\hbar}(\rho-\rho_0+\sqrt{\pi\hbar}w_1)G_{\hbar}(\rho-\rho_0+\sqrt{\pi\hbar}w_2)dw_1dw_2d\rho_0\\
&=&\int_{\mathbb{R}^{2d}}a(\rho+B(\hbar)^{-1}\rho_0)K_{\hbar}(\rho_0)d\rho_0,
\end{array}$$
where
$$K_{\hbar}(\rho_0):=\frac{1}{|\det B(\hbar)|}\int_{\mathbb{R}^{4d}}e^{-2\imath\pi\omega(w_1,w_2)}G_{\hbar}(\sqrt{\pi\hbar}w_1-B(\hbar)^{-1}\rho_0)G_{\hbar}(\sqrt{\pi\hbar}w_2-B(\hbar)^{-1}\rho_0)dw_1dw_2.$$
We start by computing $\int_{\mathbb{R}^{2d}}e^{-2\imath\pi\omega(w_1,w_2)}G(\sqrt{\pi\hbar}B(\hbar)w_1-\rho_0)dw_1$. Changing the variables, we find that it is equal to
$$\frac{e^{-2\imath\pi\left\langle Jw_2,(\pi\hbar)^{-\frac{1}{2}}B(\hbar)^{-1}\rho_0\right\rangle}}{|\det B(\hbar)|(\pi\hbar)^{d}}\int_{\mathbb{R}^{2d}}e^{-2\imath\pi\left\langle w_1,(\pi\hbar)^{-\frac{1}{2}}B(\hbar)^{-1*}Jw_2\right\rangle}G(w_1)dw_1.$$
We find then
$$K_{\hbar}(\rho_0)=\frac{2^d}{|\det B(\hbar)|}\int_{\mathbb{R}^{2d}}e^{-2\imath\pi\left\langle Jw_2,B(\hbar)^{-1}\rho_0\right\rangle}G\left(B(\hbar)^{-1*}Jw_2\right)G\left(\pi\hbar B(\hbar)w_2-\rho_0\right)dw_2$$
We make a change of variables to find that
$$K_{\hbar}(\rho_0)=2^d\int_{\mathbb{R}^{2d}}e^{-2\imath\pi\left\langle\rho_1,\rho_0\right\rangle}G\left(\rho_1\right)G\left(\pi\hbar B(\hbar)JB(\hbar)^{*}\rho_1-\rho_0\right)d\rho_1.$$
In order to verify that $a$ and $a\star\left(G_{\hbar}\sharp G_{\hbar}\right)$ have the same principal symbol, we write the Taylor formula at the point $\rho$ and find that
$$a\star\left(G_{\hbar}\sharp G_{\hbar}\right)(\rho)=a(\rho)\int_{\mathbb{R}^{2d}}K_{\hbar}(\rho_0)d\rho_0+\int_{\mathbb{R}^{2d}}K_{\hbar}(\rho_0)\int_0^1\left(d_{\rho+tB(\hbar)^{-1}\rho_0}a\right).(B(\hbar)^{-1})\rho_0dtd\rho_0.$$
We can verify that $\int_{\mathbb{R}^{2d}}K_{\hbar}(\rho_0)d\rho_0=1$. In fact, one has
$$\int_{\mathbb{R}^{2d}}\int_{\mathbb{R}^{2d}}e^{-2\imath\pi\left\langle\rho_1,\rho_0\right\rangle}G\left(\rho_1\right)G\left(\pi\hbar B(\hbar)JB(\hbar)^{*}\rho_1-\rho_0\right)d\rho_1d\rho_0=\int_{\mathbb{R}^{2d}}G(\rho_1)^2e^{2\imath\pi\langle\rho_1,A(\hbar)\rho_1\rangle}d\rho_1,$$
where $A(\hbar):=\pi\hbar B(\hbar)JB(\hbar)^{*}$. We note that $A(\hbar)$ is antisymmetric and we find that $\int_{\mathbb{R}^{2d}}K_{\hbar}(\rho_0)d\rho_0=1.$\\
Finally, we recall that we have that $\|B(\hbar)^{-1}\|_{\infty}=O(\hbar^{\gamma})$. We have to check that for a polynom $P(\rho_0)$ independent of $\hbar$, the term $\int_{\mathbb{R}^{2d}}|P|(\rho_0)|K_{\hbar}|(\rho_0)d\rho_0$ is uniformly bounded independently of $\hbar$. This quantity is bounded by
$$2^d\int_{\mathbb{R}^{2d}}\int_{\mathbb{R}^{2d}}G\left(\rho_1\right)G\left(\pi\hbar B(\hbar)JB(\hbar)^{*}\rho_1-\rho_0\right)|P|(\rho_0)|d\rho_0d\rho_1.$$
Using the fact that $\|\pi\hbar B(\hbar)JB(\hbar)^{*}\|_{\infty}$ is uniformly bounded (as $\lambda_i^+\leq\lambda_{\max}$), we have the expected property. In particular, we can use Calder\'on-Vaillancourt theorem (property~(\ref{e:calderon-vail})) to derive that $$\|\Op_{\hbar}^w(a)-\Op_{\hbar}^+(a)\|_{L^2(\mathbb{R}^d)\rightarrow L^2(\mathbb{R}^d)}=O_a(\hbar^{\gamma}).\square$$

\section{Proof of proposition~\ref{p:new-quant-torus}}\label{s:technical-lemmas}

In this appendix, we prove proposition~\ref{p:new-quant-torus} on our quantization $\Op_{\kappa}^+$. This proposition was crucial in our proof as it ensures that $\Op_{\kappa}^+$ is nonnegative and has the same nice (`product') structure as the anti-Wick quantization. We start the proof of this proposition by computing the $n$-th Fourier coefficient of $T_0(\overline{F}_1\sharp F_2)$ (where $\sharp$ is the Moyal product of two observables~\cite{DS}). We show that:
\begin{lemm} Let $F_1$ and $F_2$ be two elements in $\mathcal{S}(\mathbb{R}^{2d})$. Then, we have, for any $n$ in $\mathbb{Z}^{2d}$,
$$(T_0(\overline{F}_1\sharp F_2))_{n}:=\int_{\mathbb{T}^{2d}}e^{2\imath\pi\left\langle\rho,Jn\right\rangle}T_0(\overline{F}_1\sharp F_2)(\rho)d\rho=\left(\int_{\mathbb{R}^{2d}}e^{\frac{\imath\pi}{N}\langle n,\rho\rangle}\hat{\overline{F}}_1\left(-Jn+\rho\right)\hat{F}_2\left(-\rho\right)d\rho\right),$$
where $\hat{F}_*(\rho):=\int_{\mathbb{R}^{2d}}F_*(w)e^{-2\imath\pi\langle\rho,w\rangle}dw$ is the standard Fourier transform of $F_*$.
\end{lemm}
\emph{Proof.} Using exact expression of the Moyal product from~\cite{EZ} (see also~\cite{DS}), we write
$$T_0(\overline{F}_1\sharp F_2)(\rho)=\sum_{r\in\mathbb{Z}^{2d}}\int\int_{\mathbb{R}^{4d}}e^{-2\imath\pi\langle \rho_1,J\rho_2\rangle}\overline{F}_1\left(\frac{\rho_1}{\sqrt{2N}}+\rho+r\right)F_2\left(\frac{\rho_2}{\sqrt{2N}}+\rho+r\right)d\rho_1d\rho_2.$$
Using Poisson formula, we find that
$$T_0(\overline{F}_1\sharp F_2)(\rho)=\sum_{r\in\mathbb{Z}^{2d}}\left(\int\int\int_{\mathbb{R}^{6d}}e^{-2\imath\pi(\langle \rho_1,J\rho_2\rangle+\langle r,\rho'\rangle)}\overline{F}_1\left(\frac{\rho_1}{\sqrt{2N}}+\rho'\right)F_2\left(\frac{\rho_2}{\sqrt{2N}}+\rho'\right)d\rho_1d\rho_2d\rho'\right)e^{2\imath\pi\langle r,\rho\rangle}.$$
We recall that we are interested in the $(Jn)$-th Fourier coefficient of $T_0(\overline{F}_1\sharp F_2)$. Under the previous form, we immediatly check that
$$(T_0(\overline{F}_1\sharp F_2))_{n}=\left(\int\int\int_{\mathbb{R}^{6d}}e^{-2\imath\pi(\langle \rho_1,J\rho_2\rangle-\langle Jn,\rho'\rangle)}\overline{F}_1\left(\frac{\rho_1}{\sqrt{2N}}+\rho'\right)F_2\left(\frac{\rho_2}{\sqrt{2N}}+\rho'\right)d\rho_1d\rho_2d\rho'\right).$$
We first make the integration into the $\rho_2$ variable and we find that
$$(T_0(\overline{F}_1\sharp F_2))_{n}=\left(\int\int_{\mathbb{R}^{4d}}e^{2\imath\pi(\langle Jn,\rho'\rangle-\langle \sqrt{2N}J\rho_1,\rho'\rangle)}\overline{F}_1\left(\frac{\rho_1}{\sqrt{2N}}+\rho'\right)(2N)^d\hat{F}_2\left(-\sqrt{2N}J\rho_1\right)d\rho_1d\rho'\right).$$
Then, making the integration against the $\rho'$ variable, we find that
$$(T_0(\overline{F}_1\sharp F_2))_{n}=\left(\int_{\mathbb{R}^{2d}}e^{\frac{2\imath\pi}{\sqrt{2N}}\langle  \sqrt{2N}J\rho-Jn,\rho\rangle}\hat{\overline{F}}_1\left(-Jn+\sqrt{2N}J\rho\right)(2N)^d\hat{F}_2\left(-\sqrt{2N}J\rho\right)d\rho\right).$$
An obvious change of variables allows to find
$$(T_0(\overline{F}_1\sharp F_2))_{n}=\left(\int_{\mathbb{R}^{2d}}e^{\frac{\imath\pi}{N}\langle n,\rho\rangle}\hat{\overline{F}}_1\left(\rho-Jn\right)\hat{F}_2\left(-\rho\right)d\rho\right).\square$$

\subsection*{Proof of proposition~\ref{p:new-quant-torus}} Under the previous form, we can verify that
$$(T_0(\overline{F}_1\sharp F_2))_{n}=\sum_{r\in\mathbb{Z}^{2d}}\left(\int_{\mathbb{T}^{2d}}e^{\frac{\imath\pi}{N}\langle n,Jr\rangle}e^{\frac{\imath\pi}{N}\langle n-r,\rho\rangle}\hat{\overline{F}}_1\left(\rho-J(n-r)\right)e^{\frac{\imath\pi}{N}\langle r,\rho\rangle}\hat{F}_2\left(-\rho-Jr\right)d\rho\right).$$
We introduce the periodic function
$$\tilde{T}_{\rho}(F_{2})(\rho'):=\sum_{r\in\mathbb{Z}^{2d}}e^{\frac{\imath\pi}{N}\langle r,\rho\rangle}\hat{F}_2\left(-\rho-Jr\right)e^{-2\imath\pi\langle Jr,\rho'\rangle}.$$
Using the Poisson formula, it verifies also
\begin{equation}
\tilde{T}_{\rho}(F_{2})(\rho')=T_{\rho}(F_{2})(\rho')=\sum_{r\in\mathbb{Z}^{2d}}F_2\left(r+\rho'-\frac{J\rho}{2N}\right)e^{2\imath\pi\langle r+\rho',\rho\rangle}.
\end{equation}
With these definitions, we have
$$\overline{T_{\rho}(F_{1})}(\rho'):=\sum_{r\in\mathbb{Z}^{2d}}\overline{F}_1\left(r+\rho'-\frac{J\rho}{2N}\right)e^{-2\imath\pi\langle r+\rho',\rho\rangle}.$$
Using these new notations, we have shown the following equality which is exactly proposition~\ref{p:new-quant-torus}:
$$\Op_{\hbar}^w(T_0(\overline{F}_1\sharp F_2))=\int_{\mathbb{T}^{2d}}\Op_{\hbar}^w(T_{\rho}(F_{1}))^*\circ\Op_{\hbar}^w(T_{\rho}(F_{2}))d\rho.\square$$

\section{Proof of lemma~\ref{lemme-chat}}
\label{proof-lemma}

To complete the proof of theorem~\ref{t:main-theo}, it remains to prove lemma~\ref{lemme-chat}. To prove this lemma, we use classical properties of the entropy of a partition~\cite{Wa} (chapter $4$) that we briefly prove here (see theorem $4.3$ and $4.9$ in~\cite{Wa} for details). We fix three integers $p$, $n$ and $m$. To simplify our notations, we define the $p$-translated entropy as follows:
$$h_{2m}^p(\psi_N,\mathcal{P}):=\sum_{|\alpha|=2m}\eta\left(\mu^N(\mathbf{P}_{\alpha}^2\circ A^p)\right).$$
Mimicking the usual proof for the subadditivity of the entropy of a partition~\cite{Wa} (chapter $4$), we write
$$\begin{array}{ccc}h_{2(n+m)}^p(\psi_N,\mathcal{P})&=&-\sum_{|\alpha|=2(n+m)}\mu^N\left(\prod_{j=-m-n}^{n+m-1}P_{\alpha_j}^2\circ A^{j+p}\right)\log\mu^N\left(\prod_{j=-m+n}^{m+n-1}P_{\alpha_j}^2\circ A^{j+p}\right)\\
&+&\sum_{|\alpha|=2(n+m)}\eta\left(\frac{\mu^N\left(\prod_{j=-m-n}^{m+n-1}P_{\alpha_j}^2\circ A^{j+p}\right)}{\mu^N\left(\prod_{j=-m+n}^{m+n-1}P_{\alpha_j}^2\circ A^{j+p}\right)}\right)\mu^N\left(\prod_{j=-m+n}^{m+n-1}P_{\alpha_j}^2\circ A^{j+p}\right).
\end{array}$$
Using the concavity of the function $\eta$ and the property of partition of identity~(\ref{partition-identité-tore-classique}), we can write the following inequality:
$$h_{2(n+m)}^p(\psi_N,\mathcal{P})\leq\sum_{|\alpha|=2m}\eta\left(\mu^N\left(\prod_{j=-m+n}^{m+n-1}P_{\alpha_j}^2\circ A^{j+p}\right)\right)+\sum_{|\alpha|=2n}\eta\left(\mu^N\left(\prod_{j=-m-n}^{-m+n-1}P_{\alpha_j}^2\circ A^{j+p}\right)\right).$$
Under a more compact form, it can be reformulated as follows:
\begin{lemm} Using previous notations, one has
\begin{equation}\label{sous-additivité-chat-etape-0}
\forall p\in\mathbb{N},\ \forall n\geq 0,\ \forall m\geq 0,\ h_{2(n+m)}^p(\psi_N,\mathcal{P})\leq h_{2m}^{n+p}(\psi_N,\mathcal{P})+h_{2n}^{-m+p}(\psi_N,\mathcal{P}).
\end{equation}
\end{lemm}
We fix now two integers $m_0< m$ and write the Euclidean division $m=qm_0+r$ where $0\leq r<m_0$. We use inequality~(\ref{sous-additivité-chat-etape-0}) to derive
$$h_{2m}(\psi_N,\mathcal{P})\leq h_{2qm_0}^r(\psi_N,\mathcal{P})+h_{2r}^{-qm_0}(\psi_N,\mathcal{P}).$$
We apply one more time inequality~(\ref{sous-additivité-chat-etape-0}) to find
$$h_{2m}(\psi_N,\mathcal{P})\leq h_{2(q-1)m_0}^{r+m_0}(\psi_N,\mathcal{P})+h_{2m_0}^{-(q-1)m_0+r}(\psi_N,\mathcal{P})+h_{2r}^{-qm_0}(\psi_N,\mathcal{P}).$$
By induction, we finally have the following corollary:
\begin{coro} Using previous notations, one has
\begin{equation}\label{sous-additivité-chat-etape-finale}h_{2m}(\psi_N,\mathcal{P})\leq h_{2r}^{-qm_0}(\psi_N,\mathcal{P})+\sum_{j=1}^{q}h_{2m_0}^{-(q+1-2j)m_0+r}(\psi_N,\mathcal{P}).
\end{equation}
\end{coro}
\subsection*{Proof of lemma~\ref{lemme-chat}} This last inequality is true for any integers $(m,m_0,r)$ satisfying $m=qm_0+r$. We can now give the proof of lemma~\ref{lemme-chat}. To do this, we fix a positive integer $m_0$ and consider $(q,r)$ in $\mathbb{N}\times\mathbb{N}$ satisfying $qm_0+r=m_E(N)$ where $0\leq r<m_0$. Recall that according to Egorov property (proposition~\ref{non-exact-egorov-chat}), one has, for every $a$ in $\mathcal{C}^{\infty}(\mathbb{T}^2)$,
$$\forall\ |t|\leq m_E(N),\ \mu^N\left(a\circ A^t\right)=\mu^N(a)+o_a(1),\ \text{as}\ N\rightarrow+\infty.$$
We underline that the remainder tends to $0$ uniformly for $t$ in the allowed interval. We now apply this property to $\mathbf{P}_{\alpha}^2$ where $|\alpha|=2m_0$. Using the continuity of $\eta$, we find that
$$\forall\ |t|\leq m_E(N),\ \eta\left(\mu^N\left(\mathbf{P}_{\alpha}^2\circ A^t\right)\right)=\eta\left(\mu^N(\mathbf{P}_{\alpha}^2)\right)+o_{\alpha}(1),\ \text{as}\ N\rightarrow+\infty.$$
As $m_0$ is fixed, we can deduce from the definition of $h_{2m_0}^p(\psi_N,\mathcal{P})$ that
$$\forall\ |p|\leq m_E(N),\ h_{2m_0}^{p}(\psi_N,\mathcal{P})=h_{2m_0}(\psi_N,\mathcal{P})+o_{m_0}(1),\ \text{as}\ N\rightarrow+\infty.$$
We can apply this result in inequality~(\ref{sous-additivité-chat-etape-finale}). In this case, one has that $p=-(q+1-2j)m_0+r$ belongs to $[-m_E(N),m_E(N)]$. As $|qm_0|\leq m_E(N)$, we can also write $h_{2r}^{-qm_0}(\psi_N,\mathcal{P})=h_{2r}(\psi_N,\mathcal{P})+o_{r}(1)$ as $N$ tends to infinity. Finally, we find that
$$h_{2m_E(N)}(\psi_N,\mathcal{P})\leq h_{2r}(\psi_N,\mathcal{P})+qh_{2m_0}(\psi_N,\mathcal{P})+(q+1)R'(m_0,N),$$
where $R'(m_0,N)$ is a remainder that satisfies $\forall m_0\in\mathbb{N},\ \lim_{N\rightarrow\infty}R'(m_0,N)=0.$ The conclusion of the lemma follows from this last statement.$\square$

\end{document}